\documentclass[11pt,a4paper]{article}
\usepackage{lineno}
\pdfoutput=1

\usepackage{amsmath}
\usepackage{amssymb}
\usepackage{mathrsfs}
\usepackage{color}
\usepackage{graphicx}
\usepackage[percent]{overpic}
\usepackage{url}
\usepackage{caption}
\usepackage{subcaption}
\usepackage{enumerate}
\usepackage{overpic}
\usepackage{amsthm}
\usepackage{moreverb}
\usepackage[pdftex,colorlinks,bookmarksopen,bookmarksnumbered,citecolor=red,urlcolor=red]{hyperref}

\def\b{{\bm b}}

\def\f{{\bm f}}

\def\n{{\bm n}}
\def\u{{\bm u}}
\def\v{{\bm v}}

\def\0{\boldsymbol{0}}

\def\mubar{\overline{\mu}}
\def\ubar{\overline{\u}}
\def\vbar{\overline{\v}}

\def\qbar{\overline{q}}
\def\dt{\partial_t}

\def\cl {\nonumber \\}
\def\el {\nonumber }

\newcommand{\bm}[1]{\mbox{\boldmath{$#1$}}}
\def\div{\nabla\cdot}

\usepackage{verbatim}
\usepackage{graphicx}
\usepackage{epstopdf}
\usepackage{authblk}

\begin{document}


\date{}
\title{A novel Large Eddy Simulation model for the Quasi-Geostrophic Equations in a Finite Volume setting}
\author[1]{Michele Girfoglio\thanks{mgirfogl@sissa.it}}
\author[2]{Annalisa Quaini\thanks{quaini@math.uh.edu}}
\author[1]{Gianluigi Rozza\thanks{grozza@sissa.it}}
\affil[1]{SISSA, International School for Advanced Studies, Mathematics Area, mathLab, via Bonomea, Trieste 265 34136, Italy}
\affil[2]{Department of Mathematics, University of Houston, Houston TX 77204, USA}
\maketitle










\begin{abstract}
We present a Large Eddy Simulation (LES) approach based on a nonlinear differential low-pass 
filter for the simulation of two-dimensional barotropic flows with under-refined meshes. 
For the implementation of such model, we choose a segregated three-step algorithm combined with a
computationally efficient Finite Volume method. 
We assess the performance of our approach on the classical double-gyre wind forcing benchmark. 
The numerical experiments we present demonstrate that our nonlinear filter is an improvement over a linear filter 
since it is able to recover the four-gyre pattern of the time-averaged stream function even with
extremely coarse meshes. In addition, our LES approach provides an average kinetic energy that compares 
well with the one computed with a Direct Numerical Simulation.
\end{abstract}




\section{Introduction}

Accurate numerical simulations of geophysical flows are not only an essential tool for ocean and weather forecast,
but they could also provide insights on the mechanisms governing climate change. 
In such simulations, two-dimensional turbulence represents a major challenge.
While in the Kolmogorov three-dimensional turbulence theory \cite{Kolmogorov41-1, Kolmogorov41-2} the kinetic energy 
is transferred from the large scales to the the small scales,  in the Kraichnan-Batchelor-Leith two-dimensional turbulence 
theory \cite{Kraichnan1967, Batchelor1969, Leith1971} 
average inverse (i.e., from small to large scales) energy and direct (from large to small scales) 
enstrophy cascades are observed. 

One of the simplest models for geophysical flow is given by the Quasi-Geostrophic equations (QGE): see, e.g., \cite{Vallis2006, Cushman-Roisin2011, McWilliams2006} for mathematical and physical fundamentals, \cite{San2012, Carere2021, Strazzullo2017} for some advanced applications and \cite{QGE_review} for a recent review on this model. 
Despite the simplification in the QGE, 
when the Munk scale (a length that depends on two nondimensional quantities, the Rossby number and the
Reynolds number) is small the numerical simulation of the QGE 
becomes computationally challenging since very fine meshes are required.
In addition, often long time intervals have to be simulated, making the computational cost prohibitive. 
A way to reduce the computational cost is to adopt a Large Eddy Simulation (LES) approach that 
allows to use a coarser mesh by modeling the effect of the small scales that do not get resolved. 

In this paper, we focus on a LES model called BV-$\alpha$ \cite{Nadiga2001, Holm2003, Monteiro2015, Monteiro2014} 
that introduces a filter for the nonlinear term of the QGE model in order to correctly simulate physical flow 
when the mesh size is greater than Munk scale. Two are the main novelties of our work: i) 
the use of a nonlinear indicator function to identify the regions of the domain where the flow needs regularization
and ii) the space discretization with a computationally efficient Finite Volume method. 
Nonlinear filter stabilization has been adopted successfully for simulating two and three-dimensional incompressible flows; see, e.g., 
\cite{layton_CMAME,abigail_CMAME,BQV,Girfoglio2019}. 
However, to the best of our knowledge, it is the first time that such a methodology is utilized for 
geophysical flows. The advantage of a Finite Volume method lies in preserving
conservative quantities for the governing equations. For this reason,
Finite Volume approximations have been widely used for LES models of incompressible flows. 
Other authors have chosen to discretize the linear BV-$\alpha$ model with a Finite Difference method 
\cite{Nadiga2001, Holm2003} or a Finite Element method \cite{Monteiro2015, Monteiro2014}.

In order to assess the performance of the proposed LES approach, we consider the classical 
double-gyre wind forcing benchmark \cite{Nadiga2001, Holm2003, Monteiro2015, Monteiro2014, Greathbatch2000, San2014}. We present numerical results for test
cases: i) Rossby number $0.0036$, Reynolds number $450$ 
and ii) Rossby number $0.008$, Reynolds number  
$Re = 1000$. Most of the works on the double-gyre wind forcing benchmark consider
Rossby number $0.0016$ and Reynolds number $200$ \cite{Nadiga2001, Holm2003, Monteiro2015, Monteiro2014}, 
while flows at Rossby number $0.0036$, Reynolds number $450$ are studied in \cite{San2012,San2014}.
Our test case ii) features higher Rossby and Reynolds numbers than what studied in the literature
and is therefore more challenging.

The paper is organized as follows. In Sec.~\ref{sec:pbd}, we introduce the QGE and the BV-$\alpha$ models and the associated 
strategy for time and space discretization.
Numerical results are reported in Sec.~\ref{sec:results}, while conclusions
and future perspectives are presented in Sec.~\ref{sec:conclusions}.




\section{Problem definition}
\label{sec:pbd}

\subsection{Quasi-Geostrophic Equations}
\label{sec:NS Equations}

We consider the motion of a two-dimensional rotating 
homogeneous flow in a two-dimensional fixed domain $\Omega$ 
over a time interval of interest $(t_0, T)$. We assume that such flow can be modeled by the 
the quasi-geostrophic equations (QGE) in stream function-potential vorticity formulation, also known as the
barotropic vorticity equation (BVE).
In order to state the BVE, let $\hat{q} =  \hat{\omega} + \beta \hat{y}$ be the dimensional potential vorticity, where 
$\hat{\omega}$ is the dimensional vorticity,  $\hat{y}$ is the dimensional vertical coordinate, and $\beta$
is the gradient of the Coriolis parameter at the basin center. For convenience, we switch to 
non-dimensional variables by introducing a characteristic length $L$ and a characteristic velocity $U$. 
The non-dimensional potential vorticity $q$ is defined as:
\begin{align}
q = Ro~\omega + y, \quad  Ro = \frac{U}{\beta L^2}\label{eq:BV2}
\end{align}
where $\omega$ is the non-dimensional vorticity, $y$ is the non-dimentional vertical coordinate,
and $Ro$ is the Rossby number, which is the ratio of inertial force to Coriolis force. 
The BVE in non-dimensional variables reads:
\begin{align}
\dt q+ \div \left(\left(\nabla \times \bm{\psi}\right) q \right) - \dfrac{1}{Re} \Delta q & = F \quad \mbox{ in }\Omega \times (t_0,T), \label{eq:BV1}
\end{align}
where $\partial_t$ denotes the time derivative, $Re$ is the Reynolds number (i.e., the ratio of inertial force to viscous force), 
$F$ denotes an external forcing, and $\bm{\psi} = (0, 0, \psi)$ with  $\psi$ being the stream function. 
The kinematic relationship between vorticity $\omega$ and the streamfunction $\psi$
yields the following Poisson equation
\begin{align}
\omega = - \Delta \psi \quad \mbox{ in }\Omega \times (t_0,T). \label{eq:BV3}
\end{align}
Using \eqref{eq:BV2}, eq.~\eqref{eq:BV3} can be rewritten in terms of $q$:
\begin{align}
q = -\text{Ro} \Delta \psi + y \quad \mbox{ in }\Omega \times (t_0,T).  \label{eq:BV4}
\end{align}

To close problem \eqref{eq:BV1}, \eqref{eq:BV4}, proper boundary conditions and initial data should be provided. 
Following \cite{Nadiga2001,Holm2003,Monteiro2015,Monteiro2014,San2012}, we enforce $\psi = \omega = 0$ 
on $\partial \Omega$ and set $\omega(x,y,t_0) = 0$, which in terms of $\psi$ and $q$ become 
\begin{align}
\psi &= 0 \quad \mbox{ on }\partial \Omega \times (t_0,T), \label{eq:BV5_comp} \\
q &= y \quad \mbox{ on }\partial \Omega \times (t_0,T), \label{eq:BV5_comp2} \\
q(x,y,t_0) &= y \quad \mbox{ in }\partial \Omega. \label{eq:BV5_comp3}
\end{align}


Summarizing, the barotropic vorticity problem
is given by eqs.~\eqref{eq:BV1}, \eqref{eq:BV4} endowed with boundary conditions \eqref{eq:BV5_comp}-\eqref{eq:BV5_comp2} 
and initial data \eqref{eq:BV5_comp3}.

\subsection{The BV-$\alpha$ problem}

Despite the fact that the QGE model is a toy problem describing the main features of geophysical flows
under certain simplifying assumptions, its Direct Numerical Simulation (DNS) is still hindered by a prohibitive 
computational cost. This is especially true in the case of climate simulations that require long time intervals
(of the order of centuries). 
A DNS for the QGE model requires a mesh with mesh size smaller than the Munk scale:
\begin{align}
\delta_M = L \, \sqrt[3]{\dfrac{\text{Ro}}{\text{Re}}}. \label{eq:munk}
\end{align}
When the mesh size fails to resolve the Munk scale, the simulation provides a non-physical solution. 
A possible remedy is to introduce a model for the small (unresolved) scales in order to recover the physical solution
while resolving only the large spatial scales and so containing the computational cost.

In this paper, we propose a nonlinear variant of the so-called \emph{BV-$\alpha$ model} \cite{Nadiga2001, Holm2003, Monteiro2015, Monteiro2014} that couples the BVE model with a differential filter. Such model reads:

\begin{align}
\dt q+ \div \left(\left(\nabla \times \bm{\psi}\right) q \right) - \dfrac{1}{\text{Re}} \Delta q &= F \quad \mbox{ in }\Omega \times (t_0,T), \label{eq:BV1_comp11}\\
-\alpha^2\div \left(a(q) \nabla\overline{q}\right) +\overline{q} &= q  \quad {\rm in}~\Omega \times
(t_0,T), \label{eq:BV2_comp33} \\
-\text{Ro} \Delta \psi + y &= \overline{q} \quad \mbox{ in }\Omega \times (t_0,T), \label{eq:BV2_comp22}
\end{align}
where $\overline{q}$ is the \emph{filtered vorticity}, $\alpha$ can be interpreted as the \emph{filtering radius} 
and $a(\cdot)$ is a scalar function such that:
\begin{align*}
a(q)\simeq 0 & \mbox{ where the flow field does not need regularization;}\\
a(q)\simeq 1 & \mbox{ where the flow field does need regularization.}
\end{align*}
By setting $a(q)\equiv 1$ in  \eqref{eq:BV1_comp11}-\eqref{eq:BV2_comp22} we retrieve
the classical BV-$\alpha$ model \cite{Nadiga2001, Holm2003, Monteiro2015, Monteiro2014}. 
This model has the advantage of making the operator in the filter equations linear and constant in time, 
but we will show that its effectivity is rather limited when very coarse meshes are considered. 

Function $a$ is called \emph{indicator function} and it plays a key role in the success of the differential filter.
Taking inspiration from the large body of work on the Leray-$\alpha$ model \cite{Borggaard2009,layton_CMAME,O-hunt1988,Vreman2004,Bowers2012}, 
we propose the following indicator function:
\begin{equation}
a(q) = \dfrac{|\nabla q|}{\text{max}\left(1, ||\nabla q||_\infty\right)}. \label{eq:ind_func}
\end{equation}
Function \eqref{eq:ind_func} is mathematically convenient because of its strong monotonicity properties. 


We will refer to eq.~\eqref{eq:BV1_comp11}-\eqref{eq:BV2_comp22} with indicator function given by 
\eqref{eq:ind_func} as the nonlinear BV-$\alpha$ (or BV-$\alpha$-NL) model.

\subsection{Time and space discretization}
\label{subsec:time-discrete}

Let us start with the time discretization of problem (\ref{eq:BV1_comp11})-(\ref{eq:BV2_comp22}). 
Let $\Delta t \in \mathbb{R}$, $t^n = t_0 + n \Delta t$, with $n = 0, ..., N_T$ and $T = t_0 + N_T \Delta t$. We denote by $f^n$ the approximation of a generic quantity $f$ at the time $t^n$. 
To discretize the time derivative in \eqref{eq:BV1_comp11}
we adopt the Backward Differentiation Formula of order 1: given $q^0$, for $n \geq 0$ find the solution $(q^{n+1}, \psi^{n+1},\qbar^{n+1})$ of system: 
\begin{align}
\dfrac{1}{\Delta t} q^{n+1} + \div \left(\left(\nabla \times \bm{\psi}^{n+1}\right) q^{n+1} \right) - \dfrac{1}{Re} \Delta q^{n+1} = b^{n+1}, \label{eq:BV1_comp11_timedisc} \\ 
-\alpha^2\div \left(a^{n+1} \nabla\overline{q}^{n+1}\right) + \overline{q}^{n+1} = q^{n+1}, \label{eq:BV2_comp33_timedisc} \\
-Ro~\Delta \psi^{n+1} + y = \overline{q}^{n+1}, \label{eq:BV2_comp22_timedisc}
\end{align}
where $a^{n+1}  = a(q^{n+1} )$ and $b^{n+1} = F^{n+1} + q^n/\Delta t$. 

In order to contain the computational cost, we opt for a segregated algorithm to solve coupled problem 
(\ref{eq:BV1_comp11_timedisc})-(\ref{eq:BV2_comp22_timedisc}). A possible algorithm is as follows:
given $q^n$ and $\bm{\psi}^n$, at $t^{n+1}$ perform the following steps
 \begin{itemize}
 \item[i)] Find the vorticity $q^{n+1}$ such that
 \begin{align}
\dfrac{1}{\Delta t} q^{n+1} + \div \left(\left(\nabla \times \bm{\psi}^n \right) q^{n+1} \right) - \dfrac{1}{\text{Re}} \Delta q^{n+1} = b^{n+1}, \label{eq:BV1_comp11_timedisc_bis} 
     \end{align}
where we have replaced $\bm{\psi}^{n+1}$ in \eqref{eq:BV1_comp11_timedisc} by $\bm{\psi}^n$, i.e.~a linear extrapolation. 
 \item[ii)] Find the filtered vorticity $\overline{q}^{n+1}$ such that
  \begin{align}
 -\alpha^2\div \left(a^{n+1} \nabla\overline{q}^{n+1}\right) + \overline{q}^{n+1} = q^{n+1}. \label{eq:BV2_comp33_timedisc_bis} 
 \end{align}
 \item[iii)] Find the stream function $\psi^{n+1}$ such that
 \begin{align}
 -\text{Ro} \Delta \psi^{n+1} + y = \overline{q}^{n+1}. \label{eq:BV2_comp22_timedisc_bis}
 \end{align}
 \end{itemize}
 
For the space discretization of problem (\ref{eq:BV1_comp11_timedisc_bis})-(\ref{eq:BV2_comp22_timedisc_bis}), 
we partition the computational domain $\Omega$ into cells or control volumes $\Omega_i$,
with $i = 1, \dots, N_{c}$, where $N_{c}$ is the total number of cells in the mesh. 
We adopt a Finite Volume (FV) approximation that is derived directly from the integral form of the governing equations. 

The integral form of eq. \eqref{eq:BV1_comp11_timedisc_bis} for each volume $\Omega_i$ is given by:
\begin{align}
\frac{1}{\Delta t}\, \int_{\Omega_i} q^{n+1} d\Omega &+ \int_{\Omega_i} \div \left(\left(\nabla \times \bm{\psi}^{n}\right) q^{n+1}\right)  d\Omega \cl
&- \dfrac{1}{Re} \int_{\Omega_i} \Delta q^{n+1} d\Omega 
= \int_{\Omega_i}b^{n+1} d\Omega. \el
\end{align}
By applying the Gauss-divergence theorem, the above equation becomes:
\begin{align}\label{eq:zetFV2}
\frac{1}{\Delta t}\, \int_{\Omega_i} q^{n+1} d\Omega &+  \int_{\partial \Omega_i} \left(\left(\nabla \times \bm{\psi}^{n}\right) q^{n+1}\right) \cdot d\textbf{A} \cl
&- \dfrac{1}{Re}\int_{\partial \Omega_i} \nabla q^{n+1} \cdot d\textbf{A}  = \int_{\Omega_i}b^{n+1} d\Omega, 
\end{align}
where $\textbf{A}$ is the surface vector associated with the boundary 
of $\Omega_i$.
Then, the discretized form of eq.~\eqref{eq:zetFV2}, divided by the control volume 
$\Omega_i$, can be written as:
\begin{align}
\frac{1}{\Delta t}\, q^{n+1}_i &+ \sum_j^{} \varphi^n_j q^{n+1}_{i,j} - \dfrac{1}{Re} \sum_j^{} (\nabla q^{n+1}_i)_j \cdot \textbf{A}_j  = b^{n+1}_i, \label{eq:QGE1}
\end{align}
where $\textbf{A}_j$ is the surface vector of the $j$-th face of the control volume and
$\varphi^n_j = \left(\nabla \times \bm{\psi}_j^{n}\right) \cdot \textbf{A}_j$. 
In \eqref{eq:QGE1}, $q^{n+1}_i$ and $b^{n+1}_i$ denote the average potential vorticity
and discrete source term in control volume $\Omega_i$, while $q^{n+1}_{i,j}$
represents the potential vorticity associated to the centroid of face $j$ normalized by the volume of $\Omega_i$. 
For more details about on the treatment of the convective and diffusive terms, 
the reader is referred to \cite{Girfoglio2019, GirfoglioPSIZETA}

We deal with the space approximations of eq.~\eqref{eq:BV2_comp33_timedisc_bis} and \eqref{eq:BV2_comp22_timedisc_bis}  
in an analogous way and obtain:
\begin{align}
-\alpha^2 \sum_j a_j^{n+1} \left(\nabla\overline{q}_i^{n+1}\right)_j \cdot \textbf{A}_j + \overline{q}_i^{n+1} = q_i^{n+1}, \label{eq:QGE2}\\
  - Ro \sum_j \left(\nabla\psi_i^{n+1}\right)_j \cdot \textbf{A}_j + y_i = \overline{q}_i^{n+1} \label{eq:QGE3},
\end{align}
respectively.

In summary, the fully discretized form of problem  (\ref{eq:BV1_comp11})-(\ref{eq:BV2_comp22}) is given by system \eqref{eq:QGE1}-\eqref{eq:QGE3}.
For the implementation of the numerical scheme described in this section, we chose the FV C++ library OpenFOAM\textsuperscript{\textregistered} \cite{Weller1998}.

\section{Numerical results} \label{sec:results}

This section presents several numerical results for the QGE (i.e., no turbulence model), BV-$\alpha$, and BV-NL-$\alpha$ models. 
We consider a benchmark test that has been widely used to analyze new techniques for turbulence in geophysical flows:
the double-gyre wind forcing experiment  \cite{Nadiga2001, Holm2003, Greathbatch2000, San2011, Monteiro2015, Monteiro2014}. 

The computational domain is rectangle $[0,1] \times [-1,1]$ and the forcing is prescribed by setting $F =\sin(\pi y)$. 
We consider two different cases:  
\begin{itemize}
\item[-] Case 1: $Ro = 0.0036$ and $Re = 450$;
\item[-] Case 2: $Ro = 0.008$ and $Re = 1000$. 
\end{itemize}
The BV-$\alpha$ model has been successfully tested for $Ro = 0.0016$ and $Re = 200$ \cite{Nadiga2001, Holm2003, Monteiro2015, Monteiro2014}. 
We selected higher Rossby and Reynolds numbers because they are more challenging. We chose Case 1 because
it has been studied in \cite{San2011} with an Approximate Deconvolution model discretized with 
a Finite Element (FE) method.
However, to the best of our knowledge it is the first time that larger values of $Ro$ and $Re$ as in Case 2 are considered. 
Although the Munk scale is the same for both cases ($\delta_M/L = 0.02$), we will show that Case 2 is more critical. 
This is in line with what observed in \cite{San2011}: the simulation of a flow with higher values of $Ro$ and $Re$ 
becomes unphysical with finer meshes than the simulation of a flow with lower values of $Ro$ and $Re$ 
and equal Munk scale.
In order to validate our approach, we proceed as follows. First, we perform a QGE simulation with 
the high resolution mesh $256 \times 512$ \cite{San2011, San2014}, which has a mesh size ($h = 1/256$) 
almost 20 times smaller than the Munk scale. Then, we run experiments on two coarse meshes, 
$16 \times 32$ ($h = 1/16$) and $4 \times 8$ ($h = 1/4$), with the QGE, BV-$\alpha$, and BV-NL-$\alpha$ models and 
compare the results with the high resolution QGE solution. 
We run all the simulations from $t_0 = 0$ to $T = 100$, with time step $\Delta t = 2.5e-5$ \cite{San2014}
For all the simulations performed with the BV-$\alpha$ and BV-NL-$\alpha$ models, we  
set $\alpha = h$ following \cite{Holm2003}.

The quantities of interest for this benchmark are the time-averaged stream function $\widetilde{\psi}$
and the time-averaged potential vorticity $\widetilde{q}$ over time interval $[20,100]$, and 
the kinetic energy of the system $E$:
\begin{equation}
E = \dfrac{1}{2}\int_\Omega  \left( \left(\dfrac{\partial \psi}{\partial y}\right)^2 + \left(\dfrac{\partial \psi}{\partial x}\right)^2 \right)d\Omega.
\end{equation}\label{eq:kin_energy}
When we use the BV-NL-$\alpha$ model,
we will also compute the time-averaged indicator function $\widetilde{a}$ over time interval $[20,100]$.

\subsection{Case 1}

\begin{figure}[htb!]
\centering
\begin{subfigure}{0.193\textwidth}
         \centering
         \includegraphics[width=\textwidth]{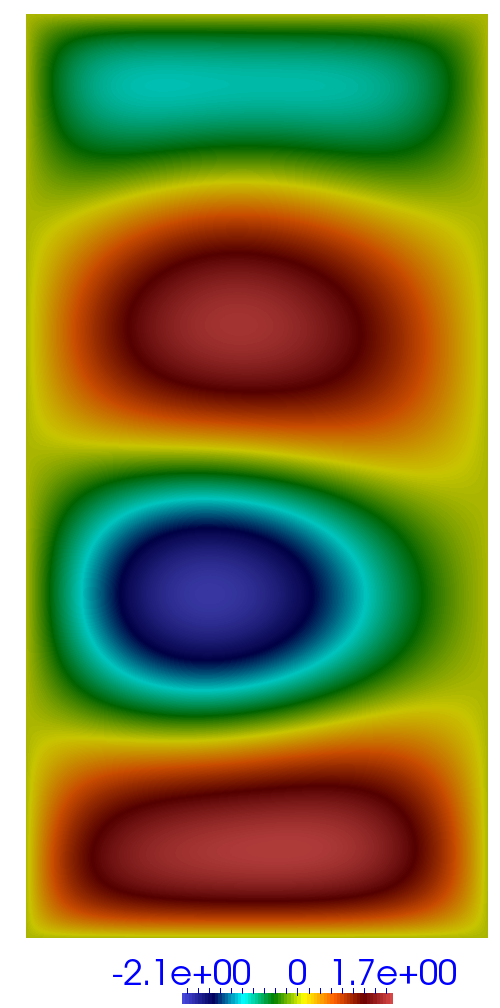}
         \caption{\scriptsize{QGE, $256\times512$}}
     \end{subfigure}
\begin{subfigure}{0.193\textwidth}
         \centering
         \includegraphics[width=\textwidth]{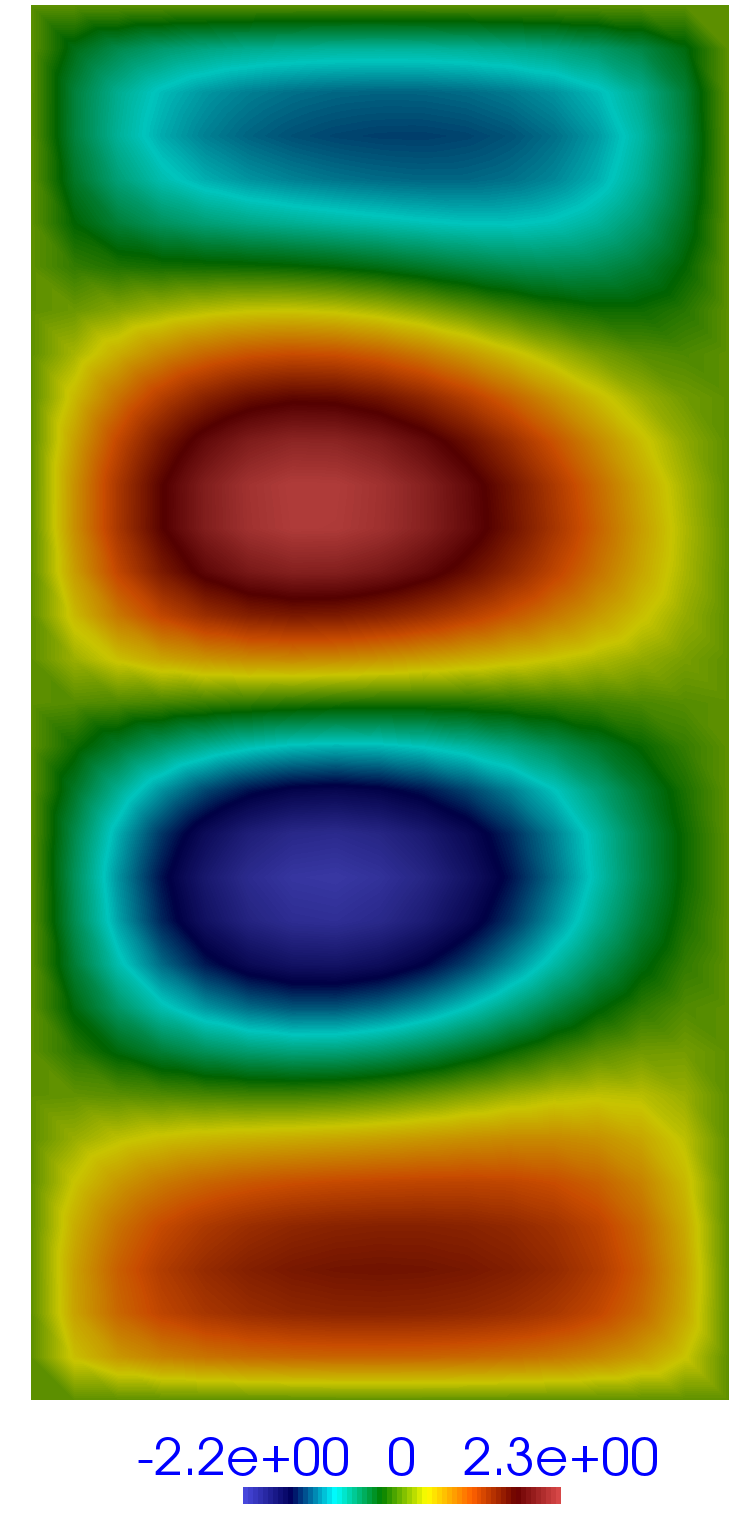}
         \caption{\scriptsize{QGE, $16 \times 32$}}
     \end{subfigure}
\begin{subfigure}{0.193\textwidth}
         \centering
         \includegraphics[width=\textwidth]{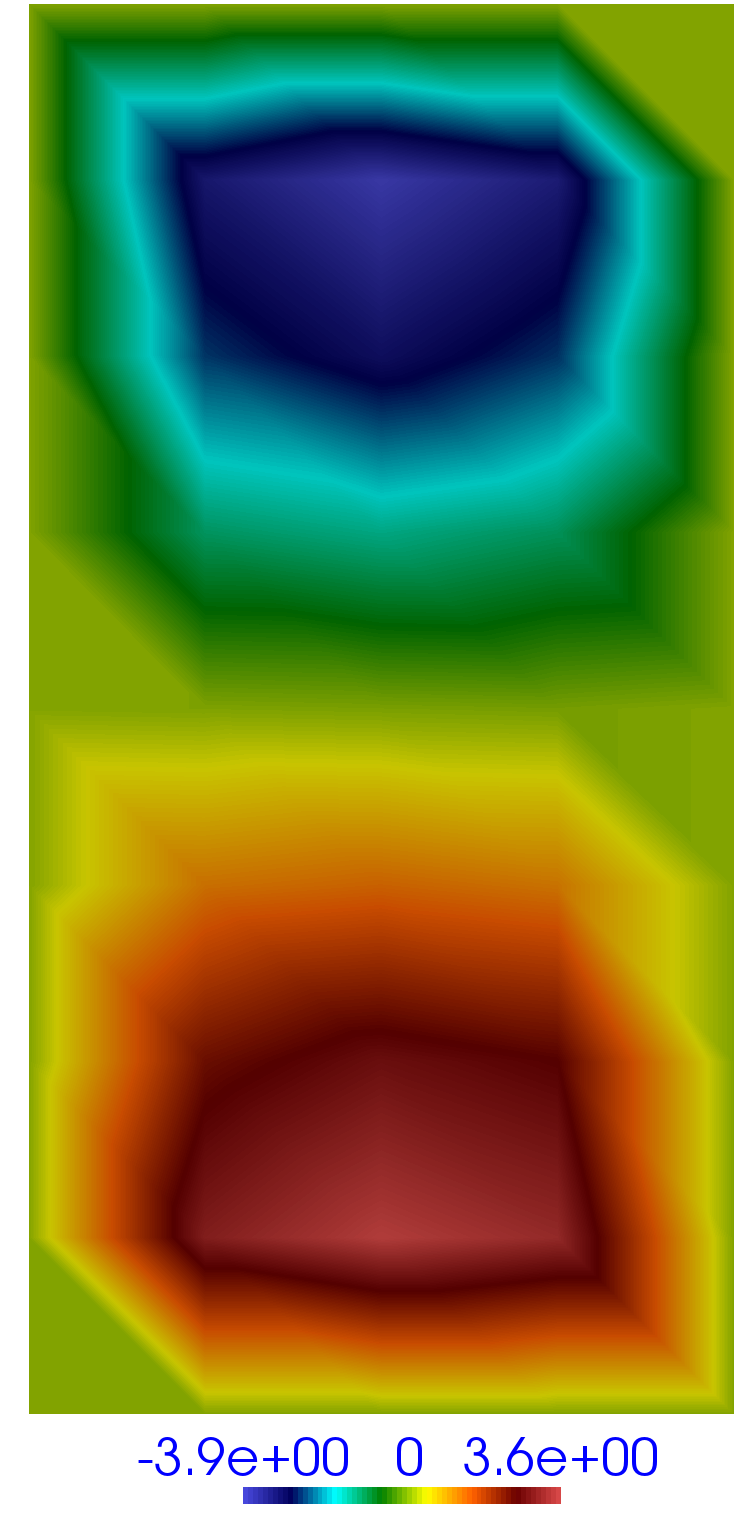}
         \caption{\scriptsize{QGE, $4 \times 8$}}
     \end{subfigure}
\begin{subfigure}{0.193\textwidth}
         \centering
         \includegraphics[width=\textwidth]{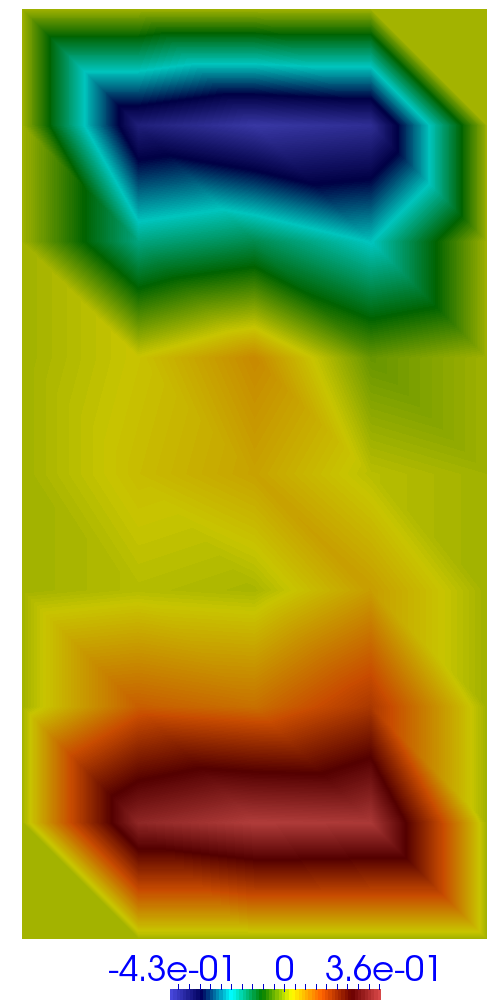}
         \caption{\scriptsize{BV-$\alpha$, $4 \times 8$}}
     \end{subfigure}
\begin{subfigure}{0.193\textwidth}
         \centering
         \includegraphics[width=\textwidth]{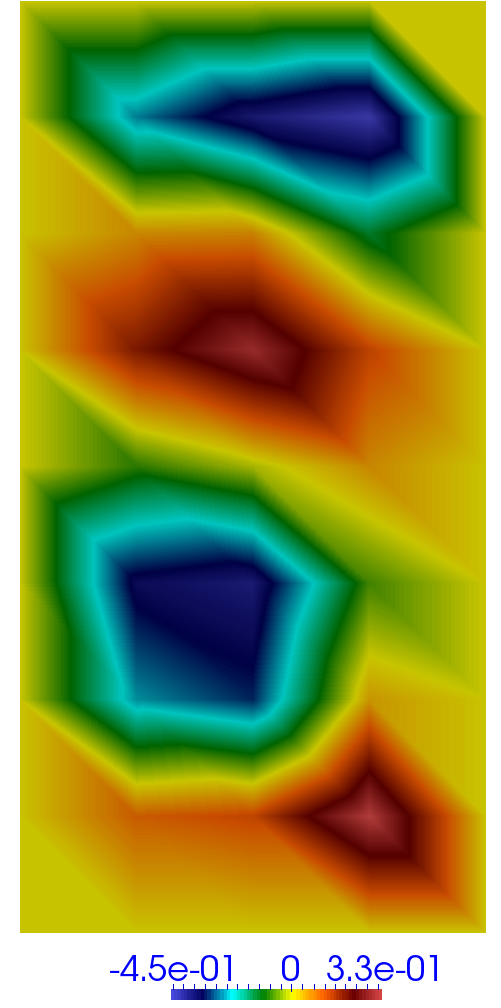}
         \caption{\scriptsize{BV-NL-$\alpha$, $4 \times 8$}}
     \end{subfigure}
\caption{Case 1: $\widetilde{\psi}$ computed with different models and different meshes. 
The specific model and mesh for each panel is reported in the corresponding subcaption.}
\label{fig:psi_first}
\end{figure}

\begin{figure}[htb!]
\centering
\begin{subfigure}{0.193\textwidth}
         \centering
         \includegraphics[width=\textwidth]{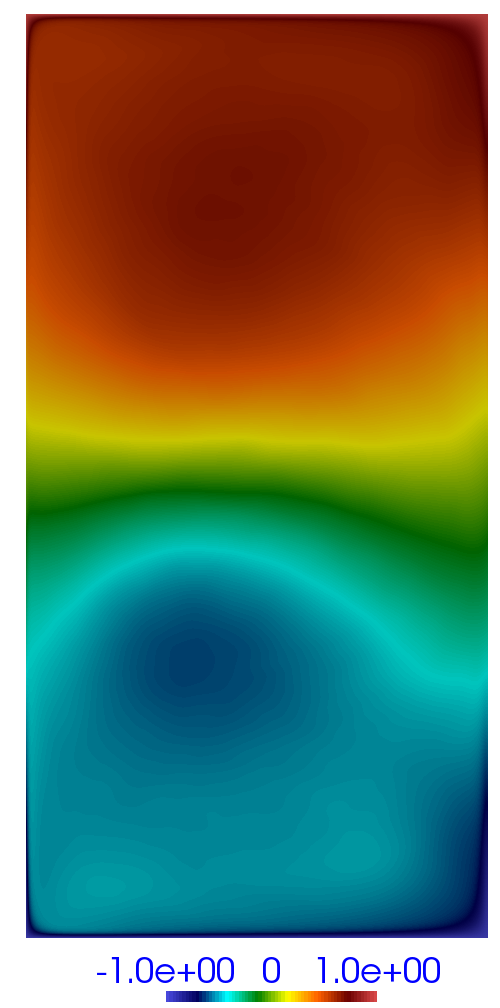}
         \caption{\scriptsize{QGE, $256\times512$}}
     \end{subfigure}
\begin{subfigure}{0.193\textwidth}
         \centering
         \includegraphics[width=\textwidth]{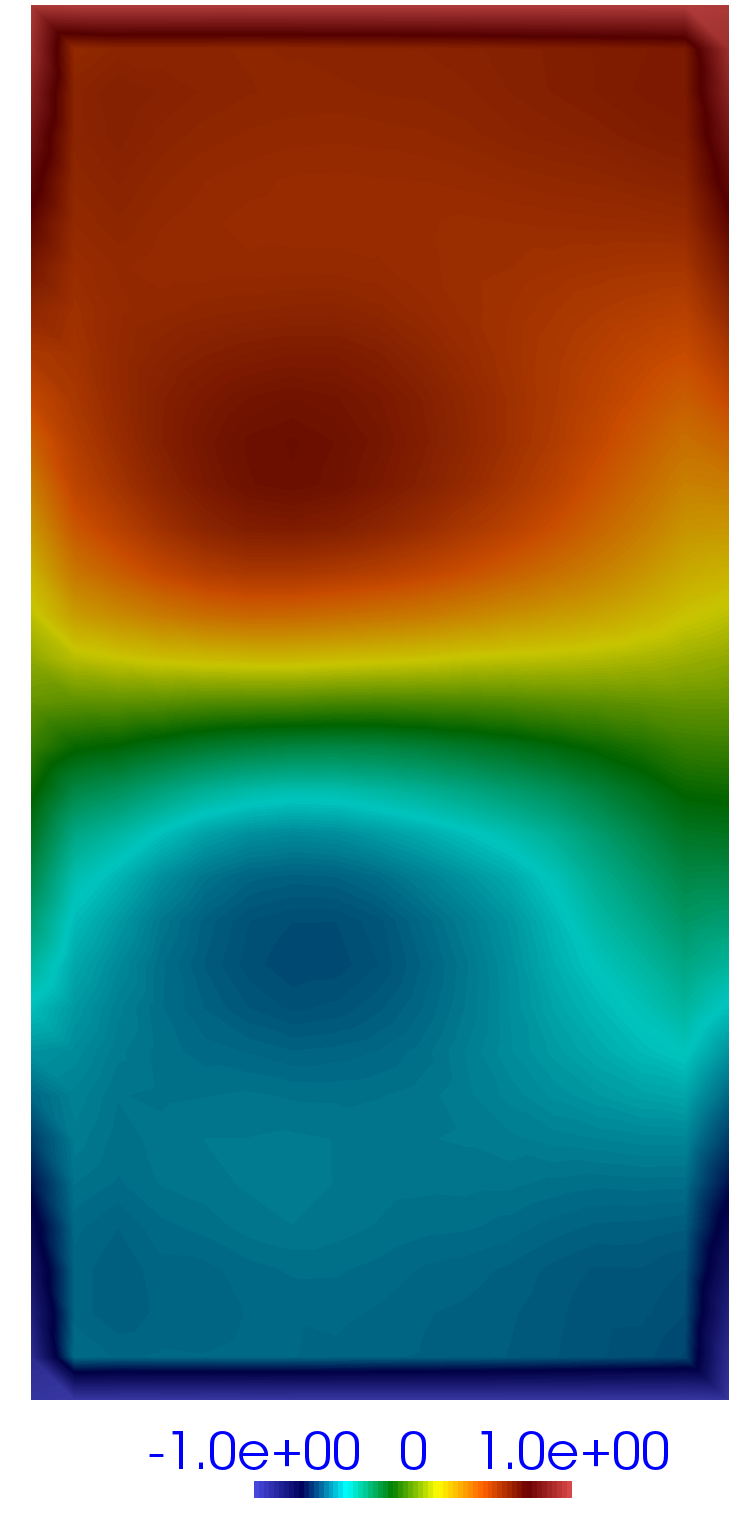}
         \caption{\scriptsize{QGE, $16 \times 32$}}
     \end{subfigure}
\begin{subfigure}{0.193\textwidth}
         \centering
         \includegraphics[width=\textwidth]{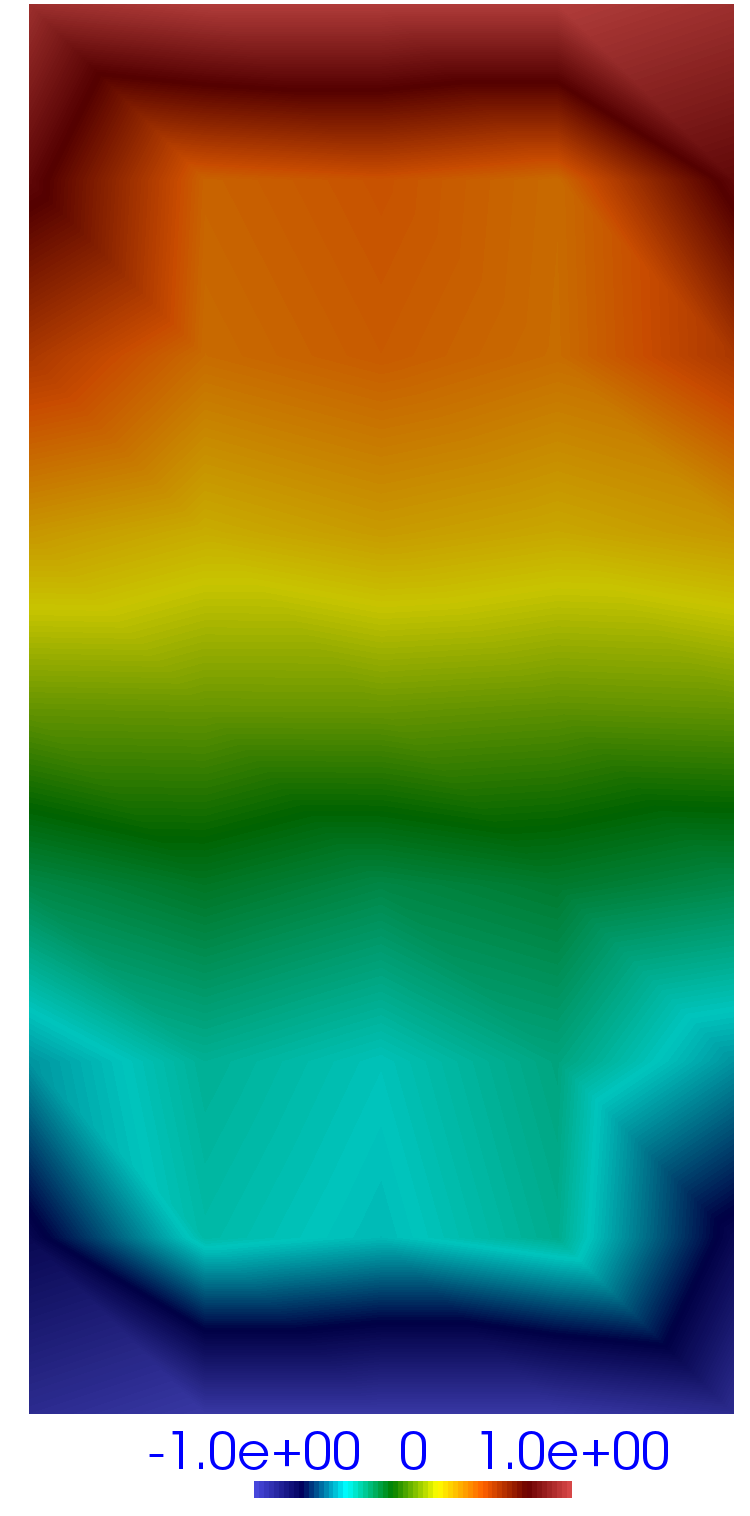}
         \caption{\scriptsize{QGE, $4 \times 8$}}
     \end{subfigure}
\begin{subfigure}{0.193\textwidth}
         \centering
         \includegraphics[width=\textwidth]{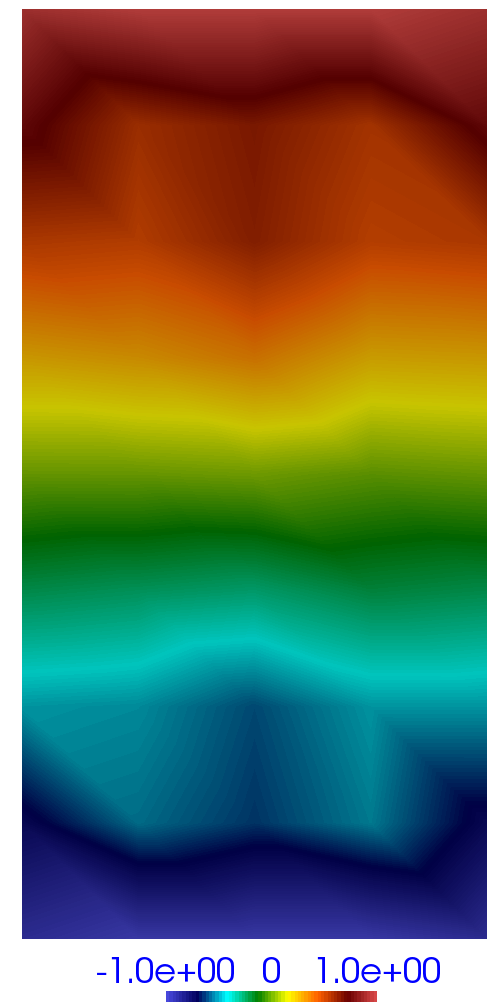}
         \caption{\scriptsize{BV-$\alpha$, $4 \times 8$}}
     \end{subfigure}
\begin{subfigure}{0.193\textwidth}
         \centering
         \includegraphics[width=\textwidth]{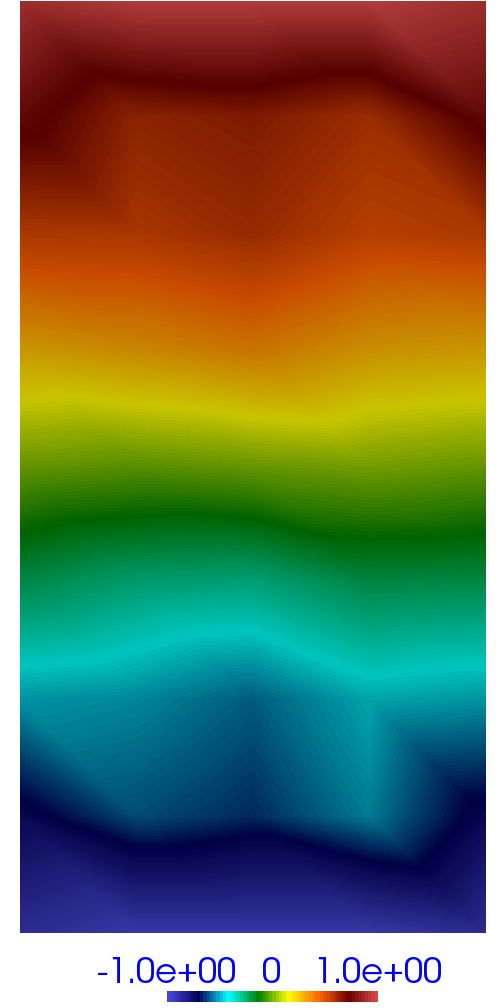}
         \caption{\scriptsize{BV-NL-$\alpha$, $4 \times 8$}}
     \end{subfigure}
\caption{Case 1:  $\widetilde{q}$ computed with different models and different meshes. 
The specific model and mesh for each panel is reported in the corresponding subcaption.
}
\label{fig:q_first}
\end{figure}

\begin{figure}[htb!]
\centering
 \begin{overpic}[width=0.48\textwidth]{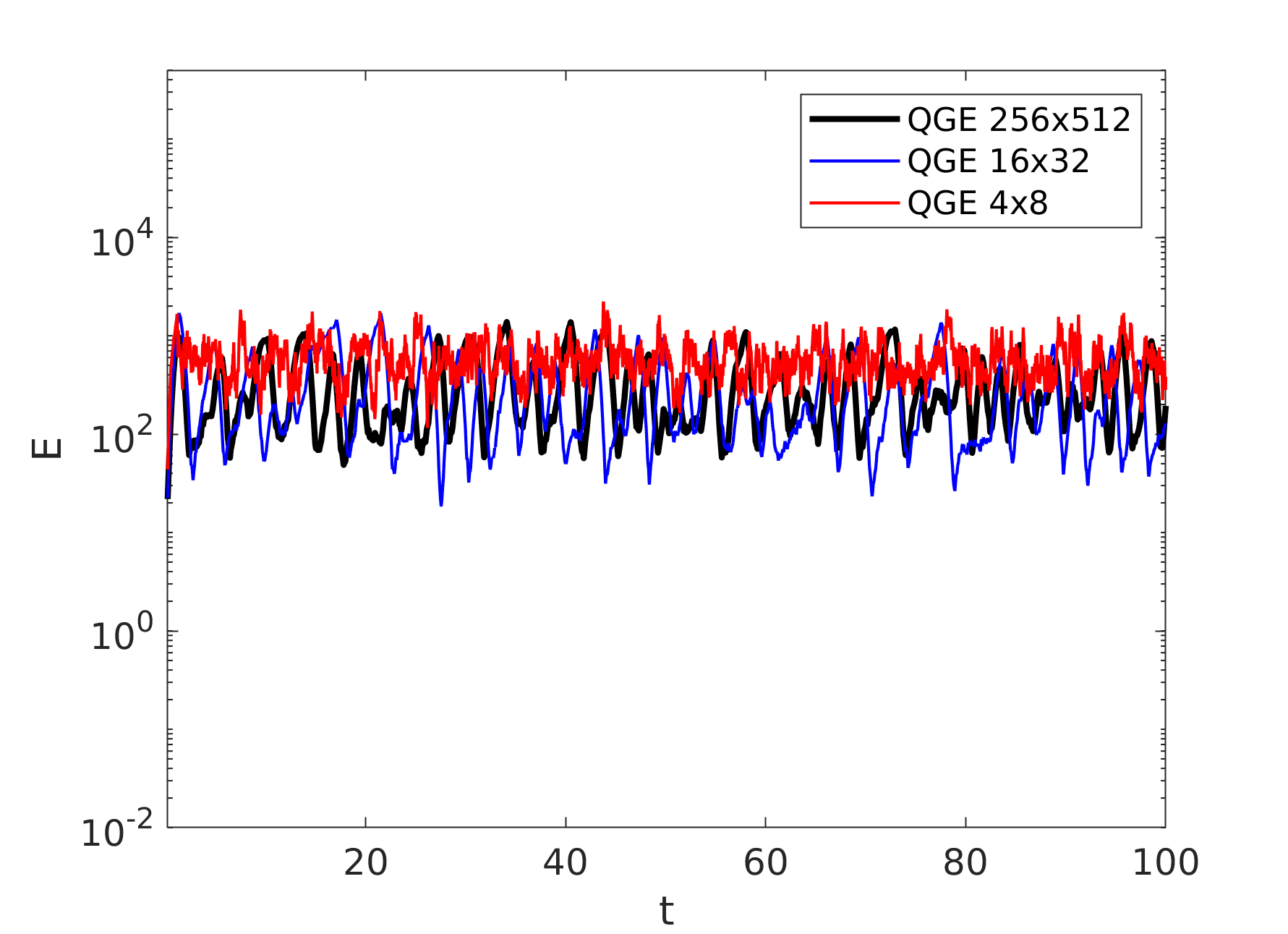}
      \end{overpic}
 \begin{overpic}[width=0.48\textwidth]{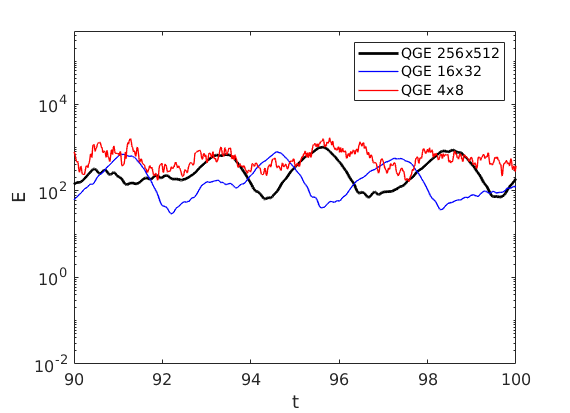}
      \end{overpic}\\
 \begin{overpic}[width=0.48\textwidth]{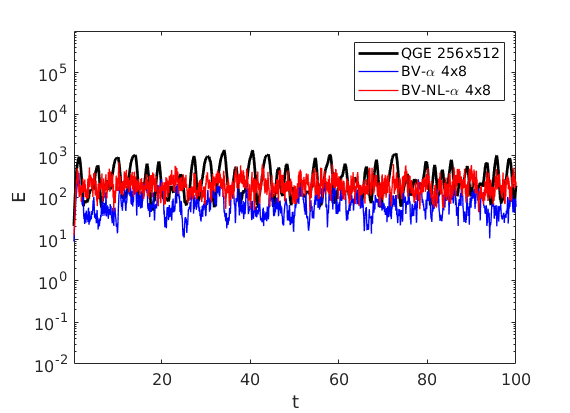}
      \end{overpic}
 \begin{overpic}[width=0.49\textwidth]{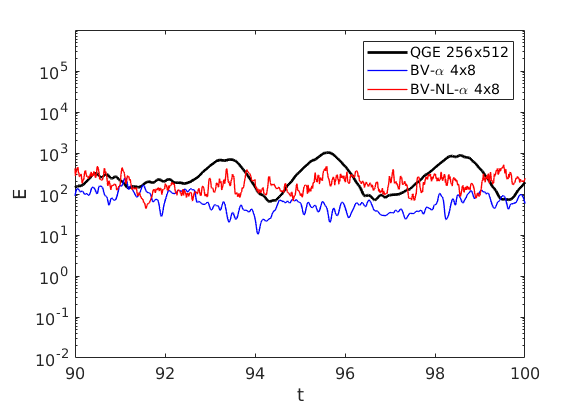}
      \end{overpic}
\caption{Case 1: time evolution of the kinetic energy \eqref{eq:kin_energy}
computed by the QGE model with all the meshes (top, left) and
corresponding zoomed-in view (top, right);  time evolution of the kinetic energy \eqref{eq:kin_energy}
computed by the QGE model on the finest, BV-$\alpha$ and BV-NL-$\alpha$ on the coarsest mesh
(bottom, left) and corresponding zoomed-in view (bottom, right).
}
\label{fig:Ek_a_first}
\end{figure}

In this section, we present the numerical results for Case 1. 
We report in  Fig.~\ref{fig:psi_first} (a) and \ref{fig:q_first} (a) the time-averaged stream function 
$\widetilde{\psi}$ and potential vorticity $\widetilde{q}$ computed with the QGE model
and high resolution mesh $256 \times 512$.
The time evolution of the kinetic energy $E$ is reported with a black line in Fig.~\ref{fig:Ek_a_first}.
We observe that our solution 
is in very good agreement with the solution in \cite{San2011, San2014}, which was computed
with a FE method. 
Hereinafter, we will refer to the solution computed with the QGE model on mesh $256 \times 512$ as the \emph{true} solution.

In \cite{San2011}, it is shown that the time-averaged stream function given by the QGE model
and computed by a FE method with the mesh $16 \times 32$ (i.e., with a mesh size about 3 times bigger than the Munk scale) 
incorrectly exhibits just two gyres instead of the four gyres seen in Fig.~\ref{fig:psi_first} (a). 
However, we observe that the QGE model approximated with our 
FV method works quite well even with mesh $16 \times 32$: it captures the four gyre pattern and the magnitude is mostly close to 
that of the true solution although some overshoots are seen. Compare Fig.~\ref{fig:psi_first} (b) and (a). 
The time-averaged potential vorticity computed with mesh $16 \times 32$ is also in good agreement
with the true $\widetilde{q}$: compare Fig.~\ref{fig:q_first} (b) and (a). For this quantity, the magnitudes
are even closer and no overshoot is observed. 
A possible reason for the better performance of our FV method with respect to a FE method could be the following:
it yields exact conservation and thus provides acceptable results despite using a
mesh with a mesh size larger than the Munk scale.

We had to push the coarseness of the mesh to $4 \times 8$ to see that the solution provided by the QGE model 
fails to show the four gyre pattern. See Fig.~\ref{fig:psi_first} (c). Notice how the maximum and minimum values of
$\widetilde{\psi}$ computed mesh $4 \times 8$ are both larger (in absolute value) than the true values. 
Similarly, the $\widetilde{q}$ computed with mesh $4 \times 8$ is not close to the true solution (see Fig.~\ref{fig:q_first} (c))
and the kinetic energy is off (see Fig.~\ref{fig:Ek_a_first}, top panels). In particular, the kinetic energy computed with mesh 
$4 \times 8$ reaches higher values for most of the time interval. 
These poor results are to be expected, since the size of mesh $4 \times 8$ is over 12 times bigger than the Munk scale. 
If we use the BV-$\alpha$ model to represent the unsolved scales with the same mesh, we observe
a slight improvement in the time-averaged potential vorticity. See Figure \ref{fig:q_first} (d). However, 
the time-averaged stream function still exhibits an incorrect pattern and its magnitude is significantly underestimated 
as shown in Figure \ref{fig:psi_first} (d). In addition, the bottom panels of Fig.~\ref{fig:Ek_a_first} show that 
the kinetic energy given by the BV-$\alpha$ model is much smaller than the true kinetic energy
over the entire time interval of interest and the frequency is off.
The BV-NL-$\alpha$ model represents an improvement over the BV-$\alpha$ model
since it is able to recover the four-gyre pattern of the time-averaged stream function (see Fig.~\ref{fig:psi_first} (e))
and it provides an average kinetic energy comparable with the true one (see Fig.~\ref{fig:Ek_a_first}, bottom panels). 
However, the magnitude of the time-averaged stream function is much smaller than it should be, as it is expected
when using a filter with such a coarse mesh, and the time-averaged vorticity is indistinguishable from the
one computed by the BV-$\alpha$ model. 


\subsection{Case 2}

\begin{figure}[htb!]
\centering
\begin{subfigure}{0.193\textwidth}
         \centering
         \includegraphics[width=\textwidth]{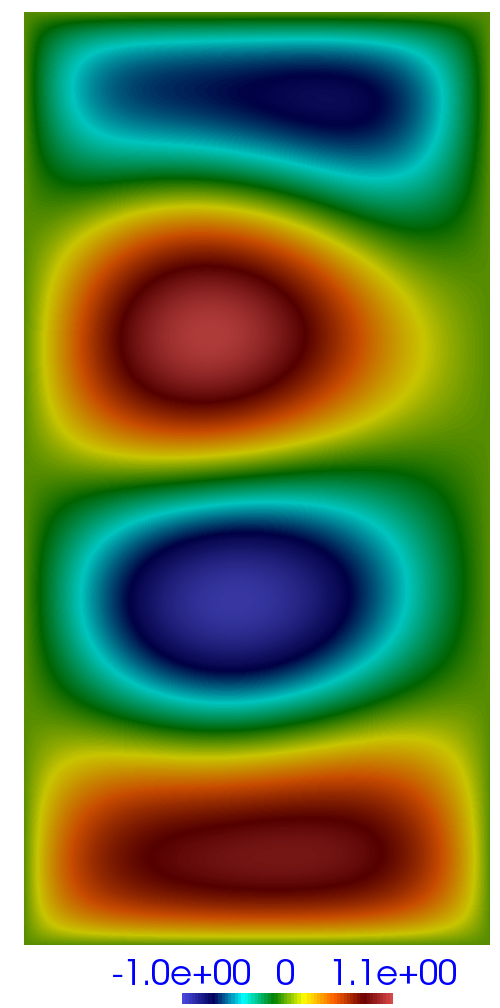}
         \caption{\scriptsize{$256\times512$}}
     \end{subfigure}
\begin{subfigure}{0.193\textwidth}
         \centering
         \includegraphics[width=\textwidth]{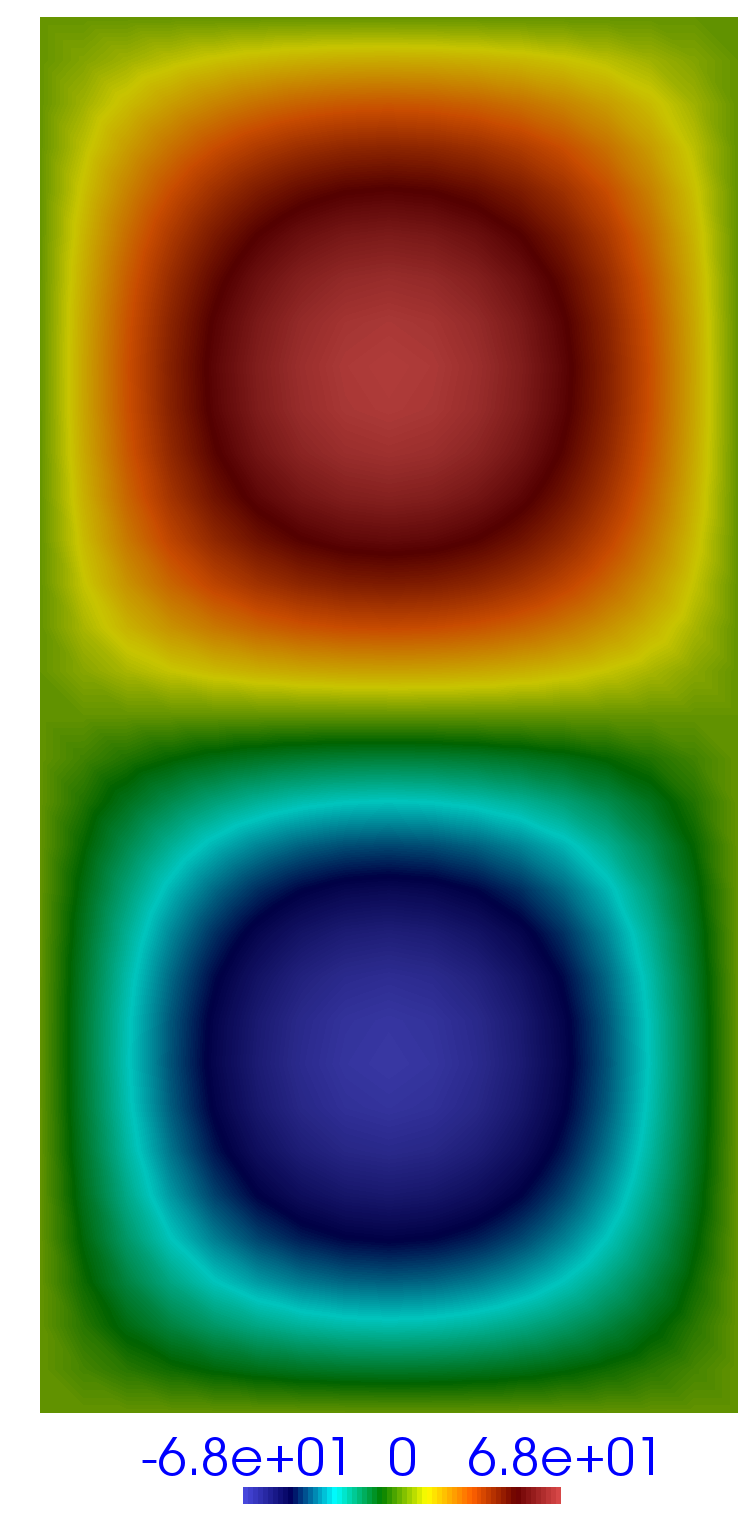}
         \caption{\scriptsize{$16 \times 32$}}
     \end{subfigure}
\begin{subfigure}{0.193\textwidth}
         \centering
         \includegraphics[width=\textwidth]{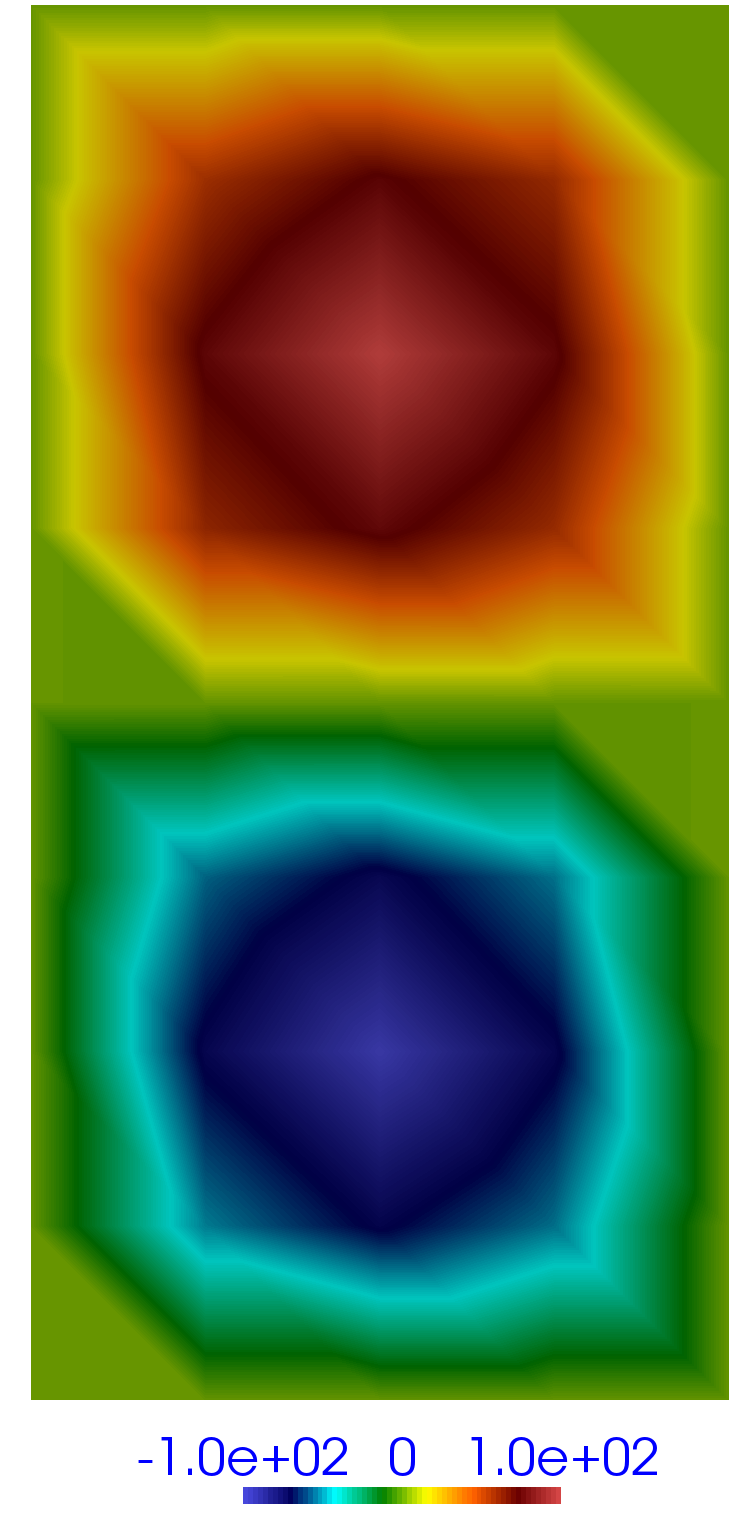}
         \caption{\scriptsize{$4 \times 8$}}
     \end{subfigure}
\caption{Case 2:  $\widetilde{\psi}$ computed by the QGE model with different meshes.}
\label{fig:psi_second_QGE}
\end{figure}

\begin{figure}[htb!]
\centering
\begin{subfigure}{0.193\textwidth}
         \centering
         \includegraphics[width=\textwidth]{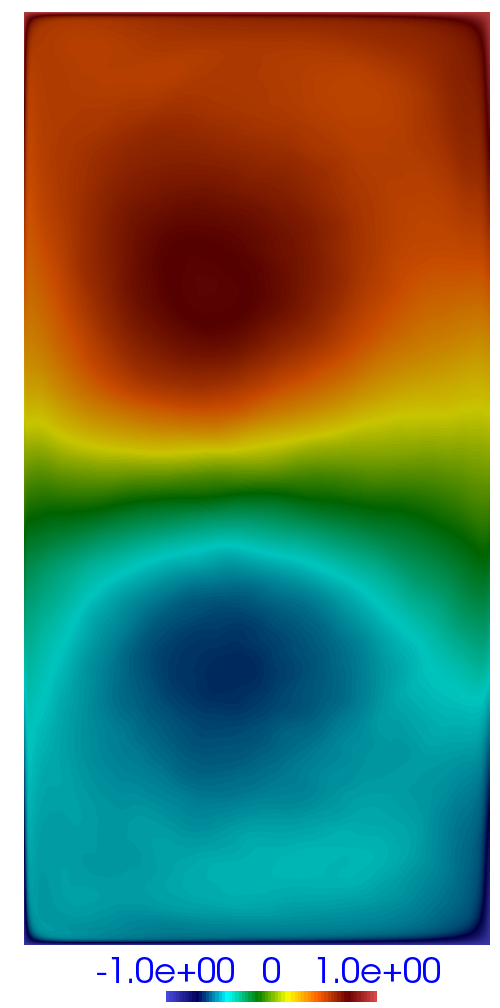}
         \caption{\scriptsize{$256\times512$}}
     \end{subfigure}
\begin{subfigure}{0.193\textwidth}
         \centering
         \includegraphics[width=\textwidth]{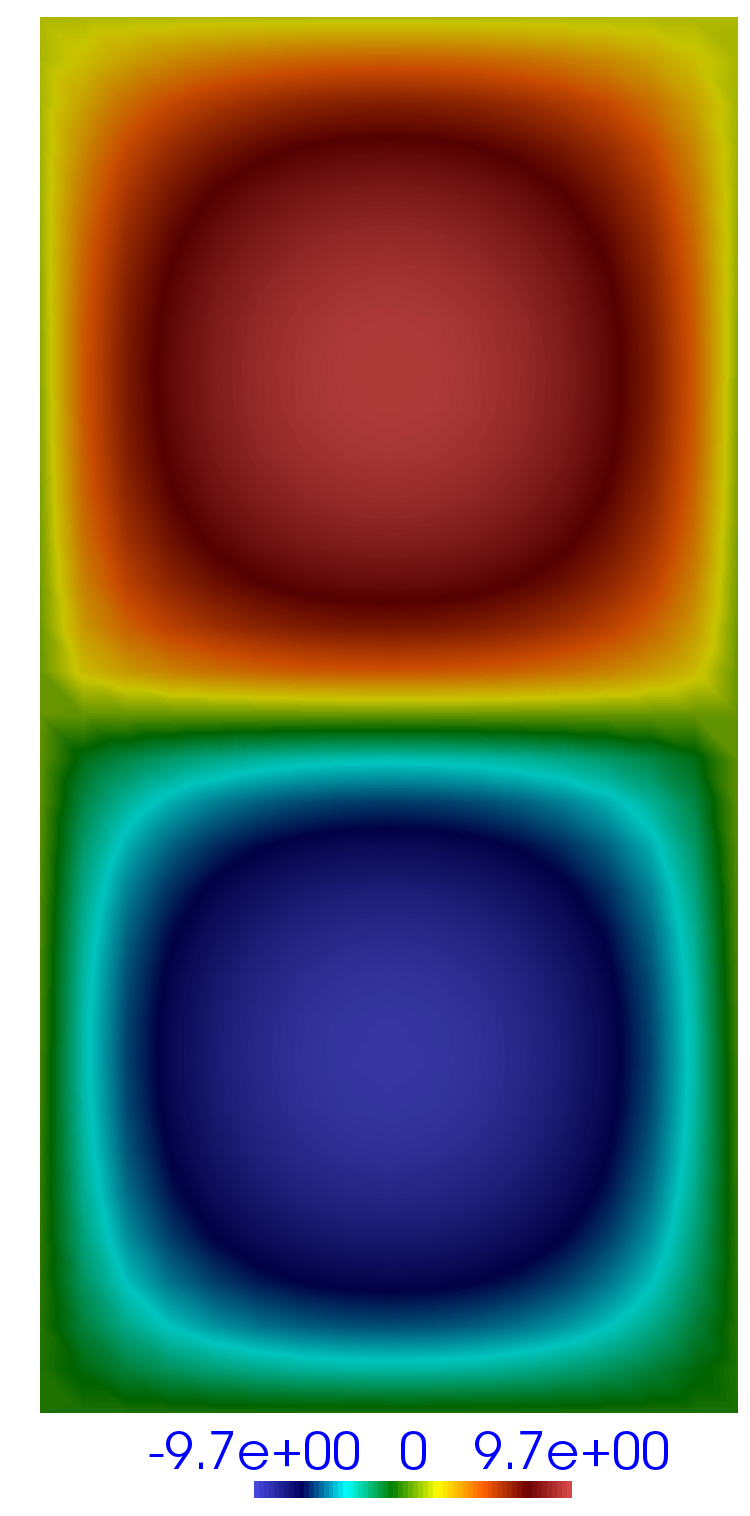}
         \caption{\scriptsize{$16 \times 32$}}
     \end{subfigure}
\begin{subfigure}{0.193\textwidth}
         \centering
         \includegraphics[width=\textwidth]{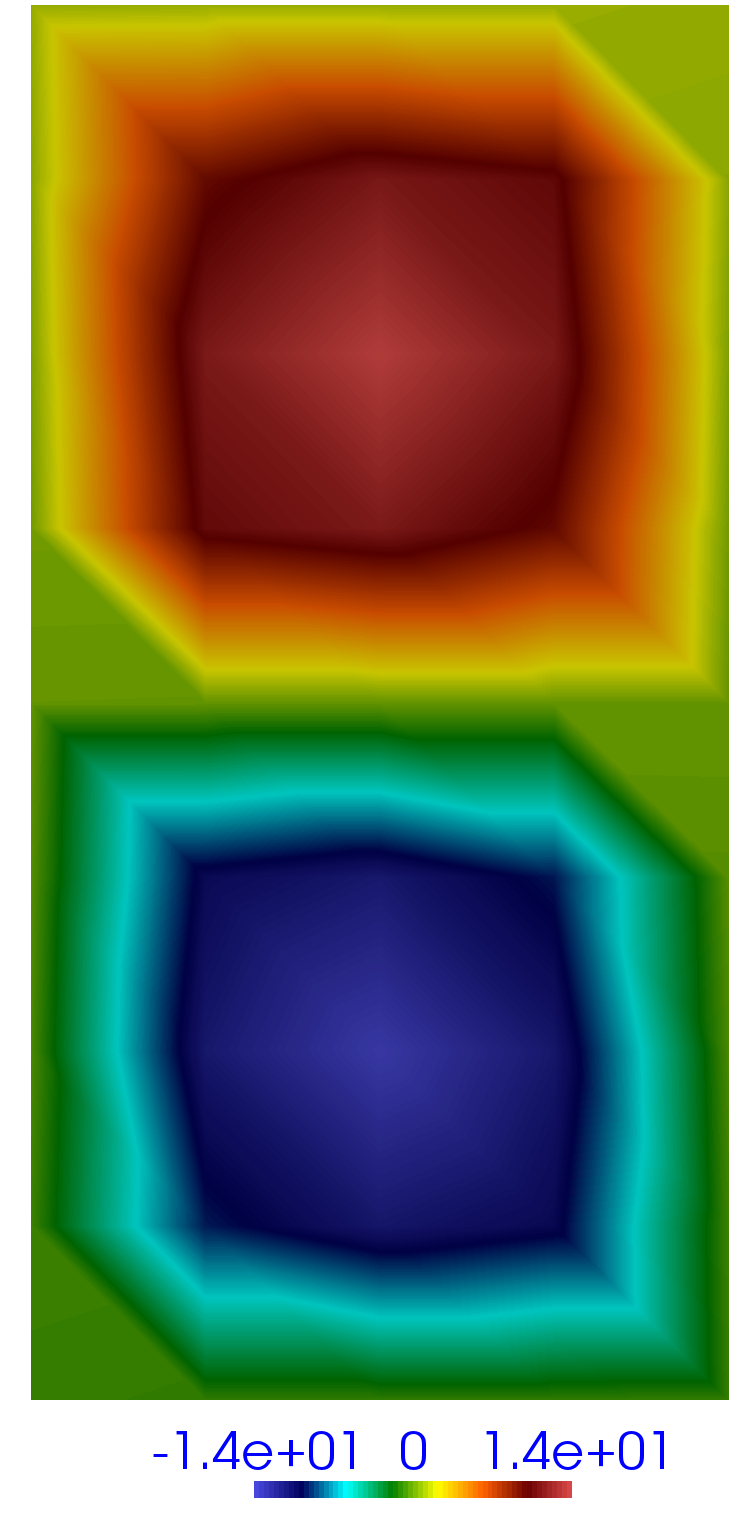}
         \caption{\scriptsize{$4 \times 8$}}
     \end{subfigure}
\caption{Case 2:  $\widetilde{q}$ computed by the QGE model with different meshes.}
\label{fig:q_second_QGE}
\end{figure}

\begin{figure}
\centering
 \begin{overpic}[width=0.48\textwidth]{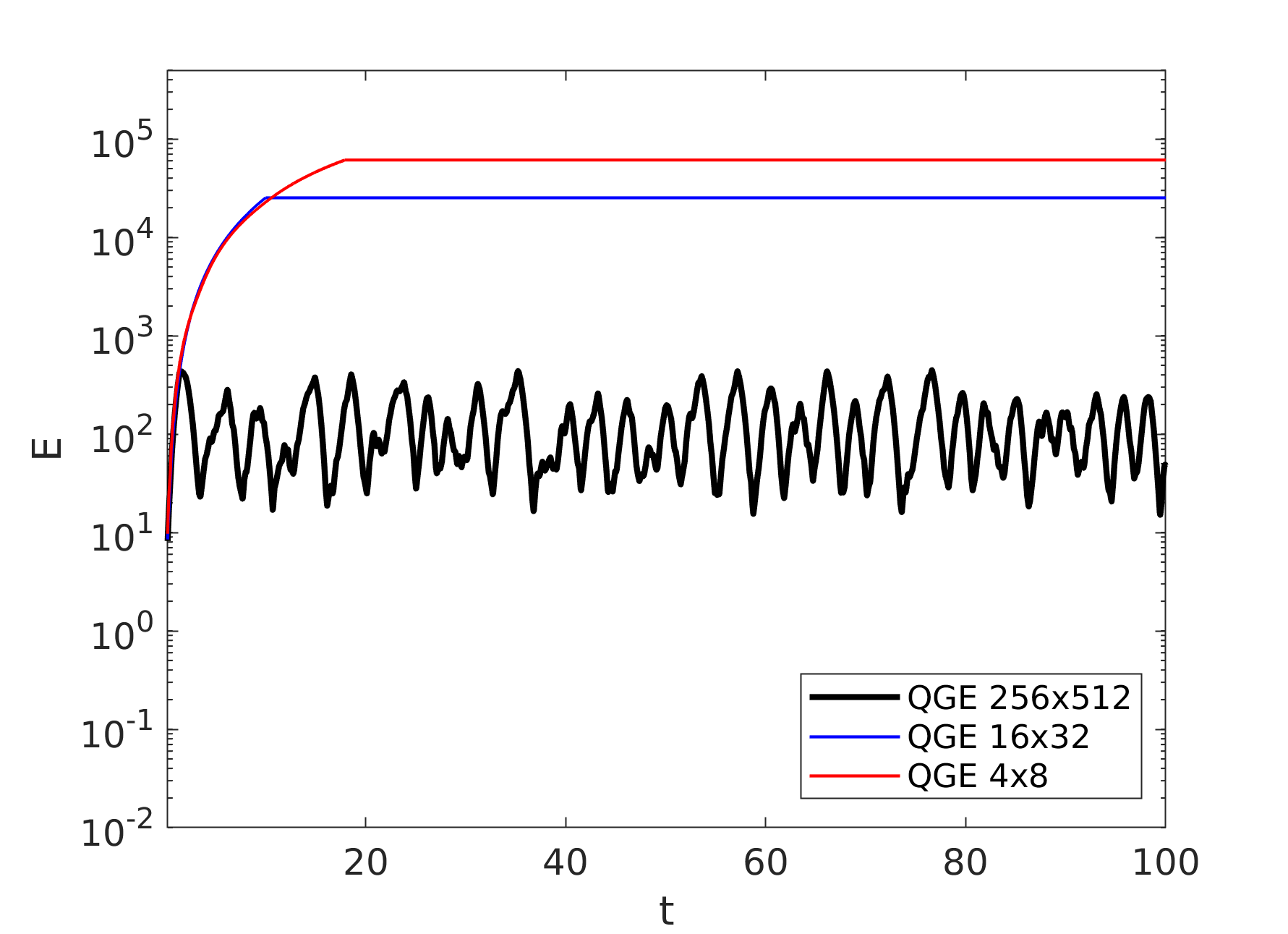}
      \end{overpic}
\caption{Case 2: time evolution of the kinetic energy \eqref{eq:kin_energy}
computed by the QGE model with different meshes.}
\label{fig:Ek_a_second_QGE}
\end{figure}

Let us now turn to Case 2.
Once again we start with the QGE model and mesh $256 \times 512$. 
Fig.~\ref{fig:psi_second_QGE} (a) and \ref{fig:q_second_QGE} (a)
display the time-averaged stream function $\widetilde{\psi}$ and potential vorticity $\widetilde{q}$, 
respectively. We will refer to this solution as the \emph{true} solution.
The corresponding  kinetic energy $E$ \eqref{eq:kin_energy} is shown in Fig.~\ref{fig:Ek_a_second_QGE} 
with a black line.

Next, we stick to the QGE model but consider coarser meshes. 
We notice that while in Case 1 the QGE model with mesh $16 \times 32$
provided a physical solution, that is it no more true for Case 2. In fact, by comparing
Fig.~\ref{fig:psi_second_QGE} (b) with Fig.~\ref{fig:psi_second_QGE} (a) 
we see that the time-averaged stream function computed with mesh $16 \times 32$
exhibits an incorrect two gyre pattern and a magnitude almost 70 times larger than the true solution. 
Similarly, the time-averaged potential vorticity shows a different pattern and a much larger magnitude.
Compare Fig.~\ref{fig:q_second_QGE} (b) and (a). In addition, from Fig.~\ref{fig:Ek_a_second_QGE}
we see that the kinetic energy computed with mesh $16 \times 32$ is not oscillatory 
and takes values much larger (up to about two orders of magnitude) than the true kinetic energy.
Poorer results are obtained with mesh $4 \times 8$. 
See Fig.~\ref{fig:psi_second_QGE} (c), \ref{fig:q_second_QGE} (c), and \ref{fig:Ek_a_second_QGE}.

\begin{figure}[htb!]
\centering
\begin{subfigure}{0.193\textwidth}
         \centering
         \includegraphics[width=\textwidth]{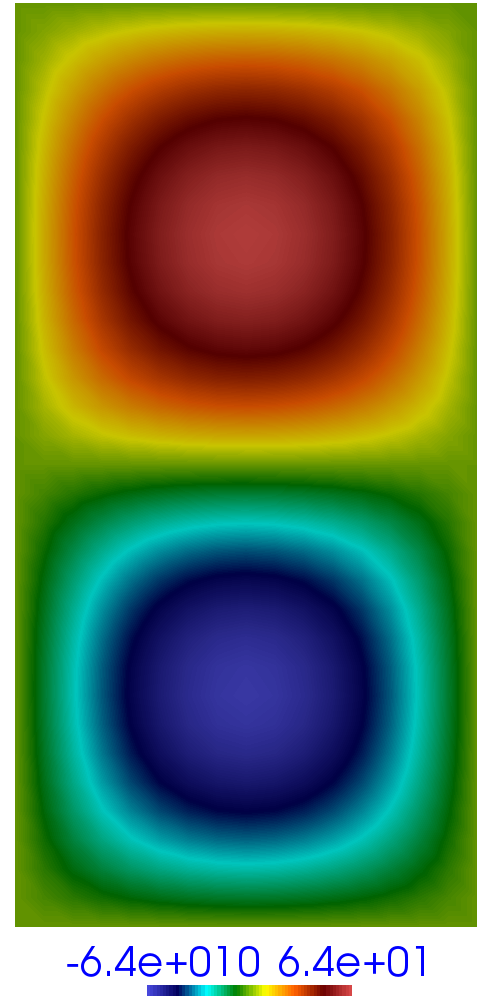}
         \caption{\scriptsize{$\widetilde{\psi}$, $16 \times 32$}}
     \end{subfigure}
\begin{subfigure}{0.193\textwidth}
         \centering
         \includegraphics[width=\textwidth]{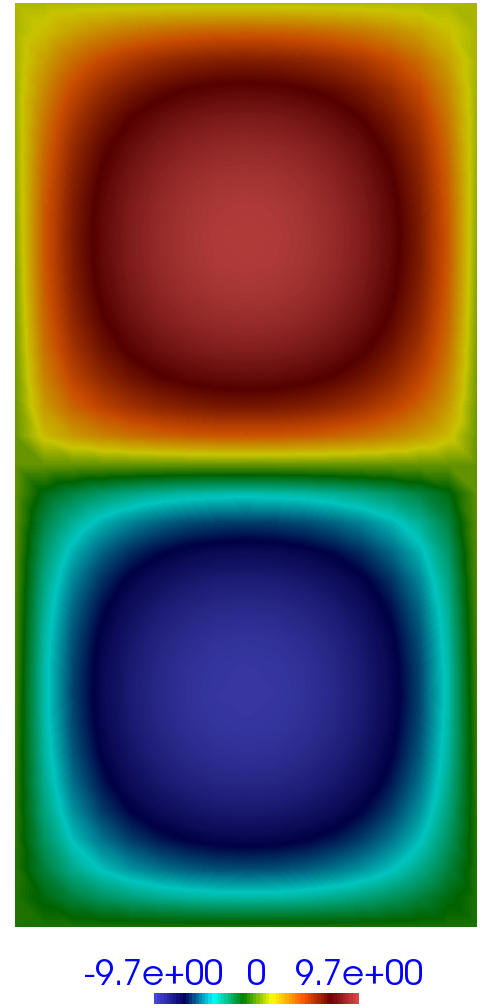}
         \caption{\scriptsize{$\widetilde{q}$, $16 \times 32$}}
     \end{subfigure}
\begin{subfigure}{0.193\textwidth}
         \centering
         \includegraphics[width=\textwidth]{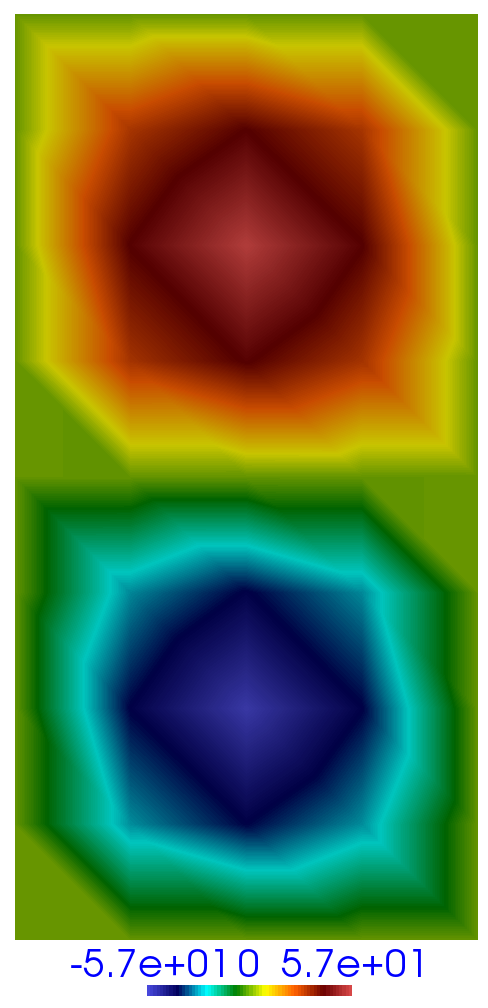}
         \caption{\scriptsize{$\widetilde{\psi}$, $4 \times 8$}}
     \end{subfigure}
     \begin{subfigure}{0.193\textwidth}
         \centering
         \includegraphics[width=\textwidth]{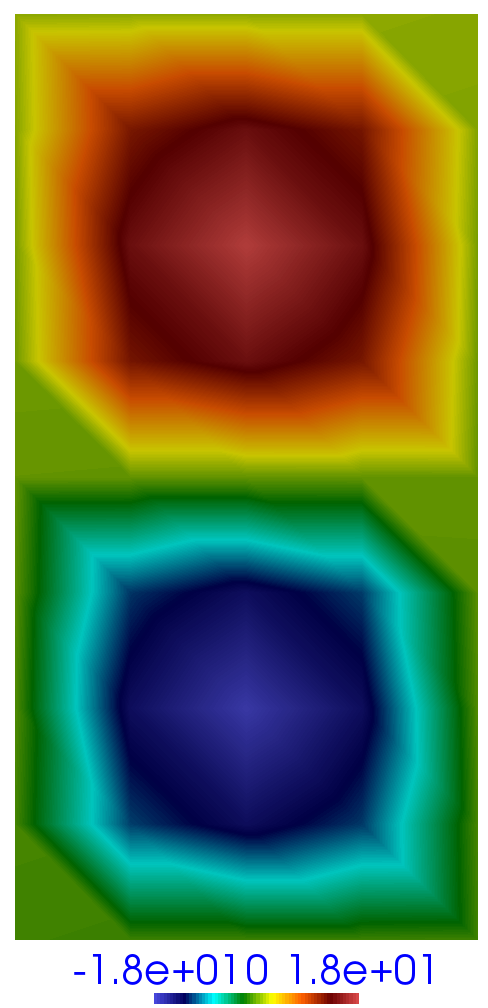}
         \caption{\scriptsize{$\widetilde{q}$, $4 \times 8$}}
     \end{subfigure}
\caption{Case 2: $\widetilde{\psi}$ and $\widetilde{q}$ computed with the BV-$\alpha$ model
and meshes $16 \times 32$ and  $4 \times 8$. }
\label{fig:BV_second_filter}
\end{figure}

\begin{figure}[htb!]
\centering
       \begin{overpic}[width=0.48\textwidth]{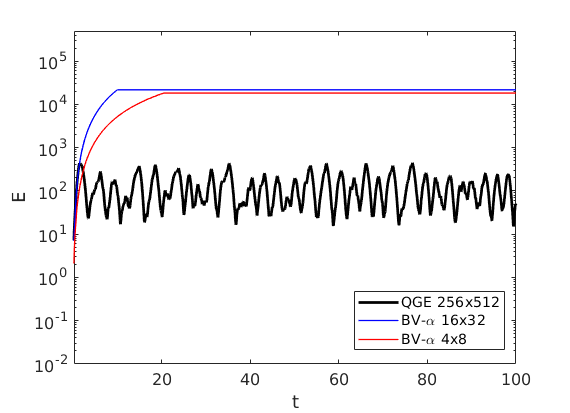}
      \end{overpic}
\caption{Case 2: time evolution of the kinetic energy \eqref{eq:kin_energy}
computed by the BV-$\alpha$ model with meshes $16 \times 32$ and  $4 \times 8$ and compared to the
true kinetic energy.} 
\label{fig:Ek_a_second_BVa}
\end{figure}

\begin{figure}[htb!]
\centering
\begin{subfigure}{0.193\textwidth}
         \centering
         \includegraphics[width=\textwidth]{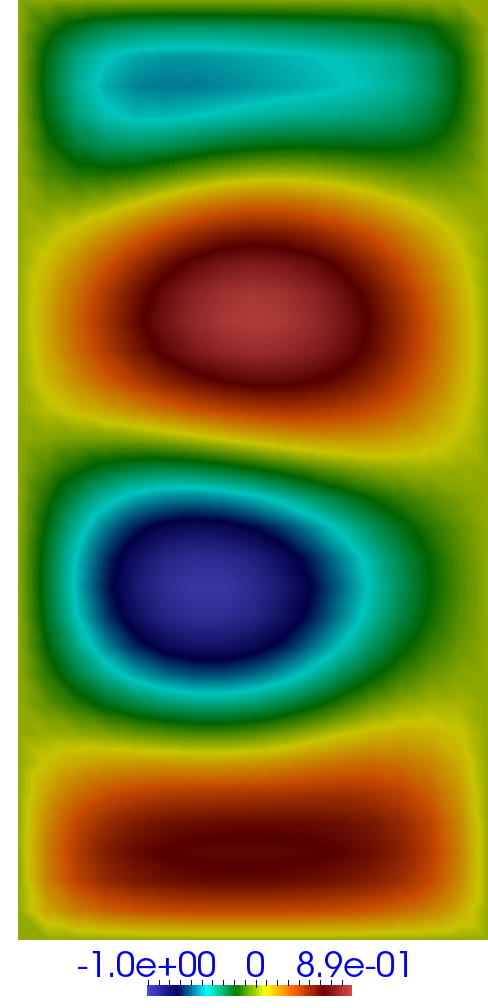}
         \caption{\scriptsize{$\widetilde{\psi}$, $16 \times 32$}}
     \end{subfigure}
\begin{subfigure}{0.193\textwidth}
         \centering
         \includegraphics[width=\textwidth]{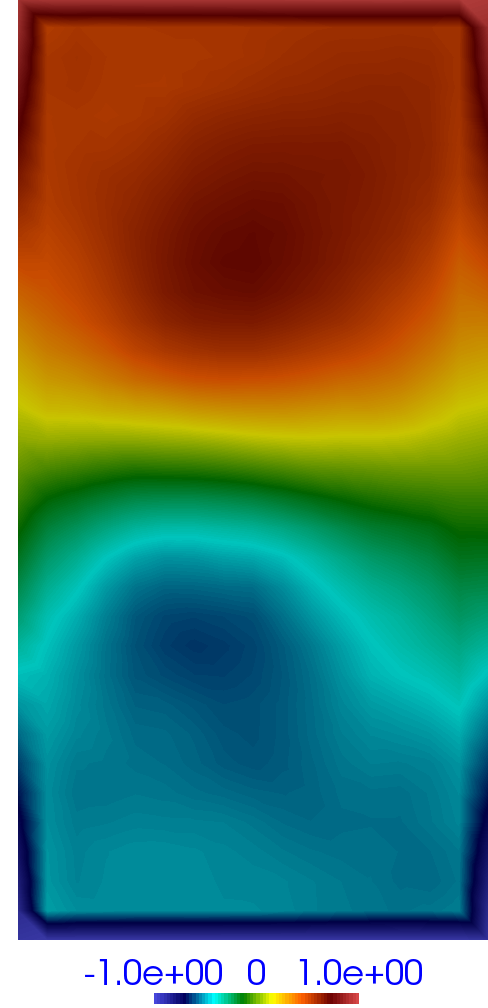}
         \caption{\scriptsize{$\widetilde{q}$, $16 \times 32$}}
     \end{subfigure}
\begin{subfigure}{0.193\textwidth}
         \centering
         \includegraphics[width=\textwidth]{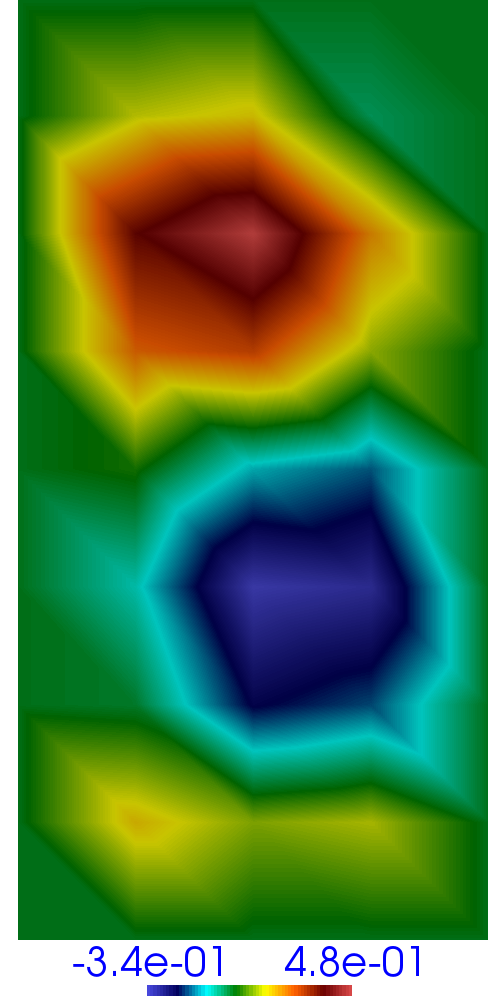}
         \caption{\scriptsize{$\widetilde{\psi}$, $4 \times 8$}}
     \end{subfigure}
     \begin{subfigure}{0.193\textwidth}
         \centering
         \includegraphics[width=\textwidth]{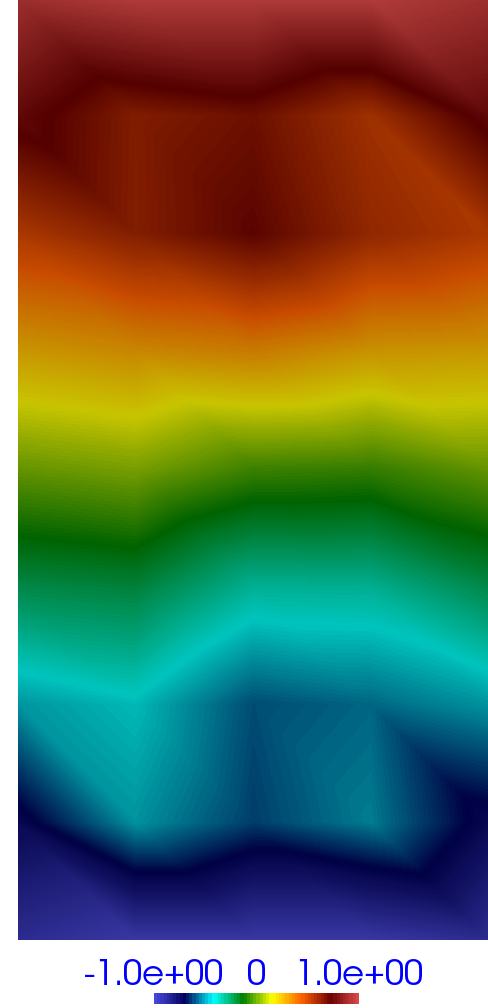}
         \caption{\scriptsize{$\widetilde{q}$, $4 \times 8$}}
     \end{subfigure}
\caption{Case 2: $\widetilde{\psi}$ and $\widetilde{q}$ computed with the BV-NL-$\alpha$ model
and meshes $16 \times 32$ and  $4 \times 8$. }
\label{fig:BV_NL_second_filter}
\end{figure}

\begin{figure}
\centering
       \begin{overpic}[width=0.48\textwidth]{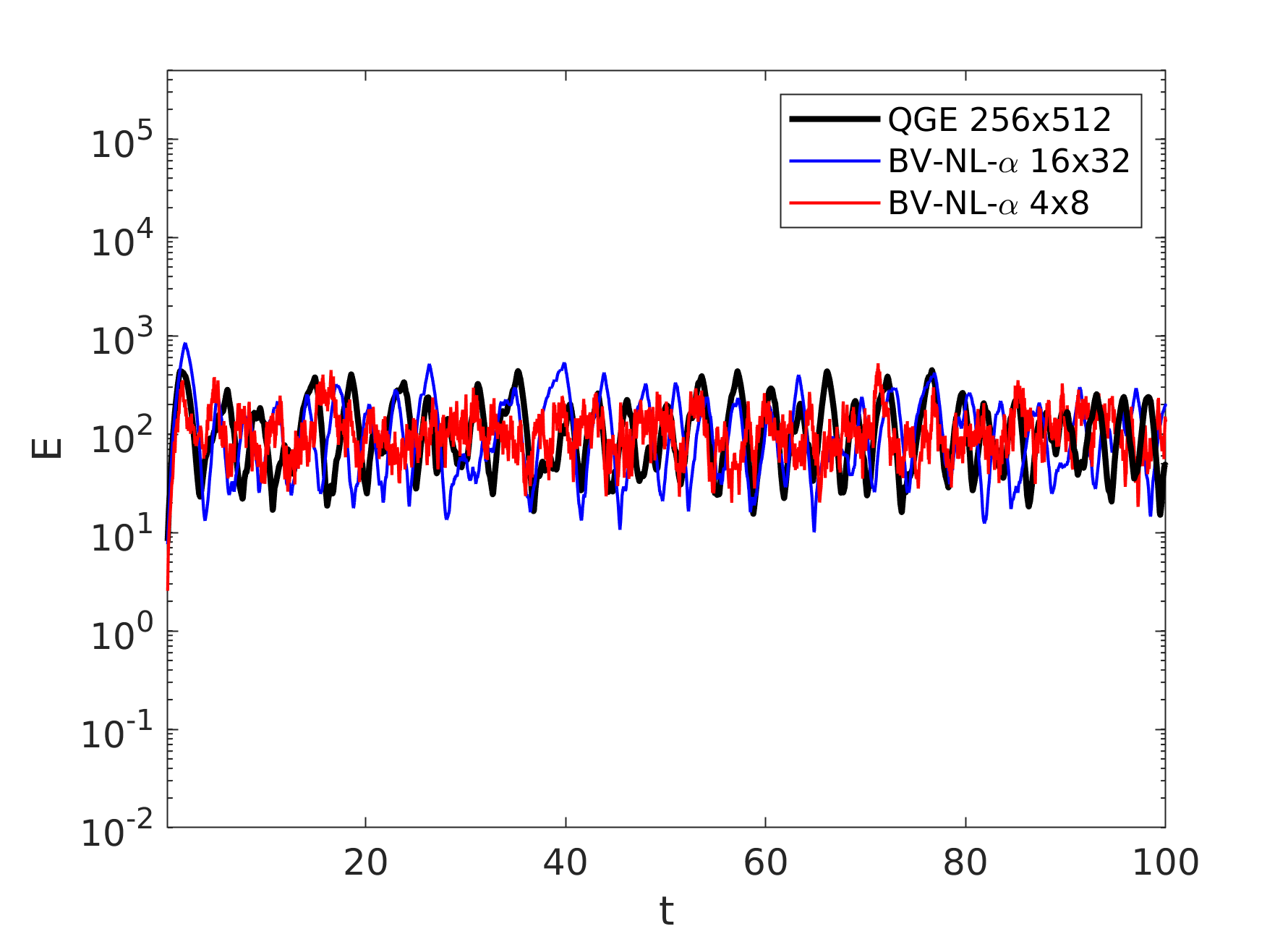}
      \end{overpic}
\begin{overpic}[width=0.49\textwidth]{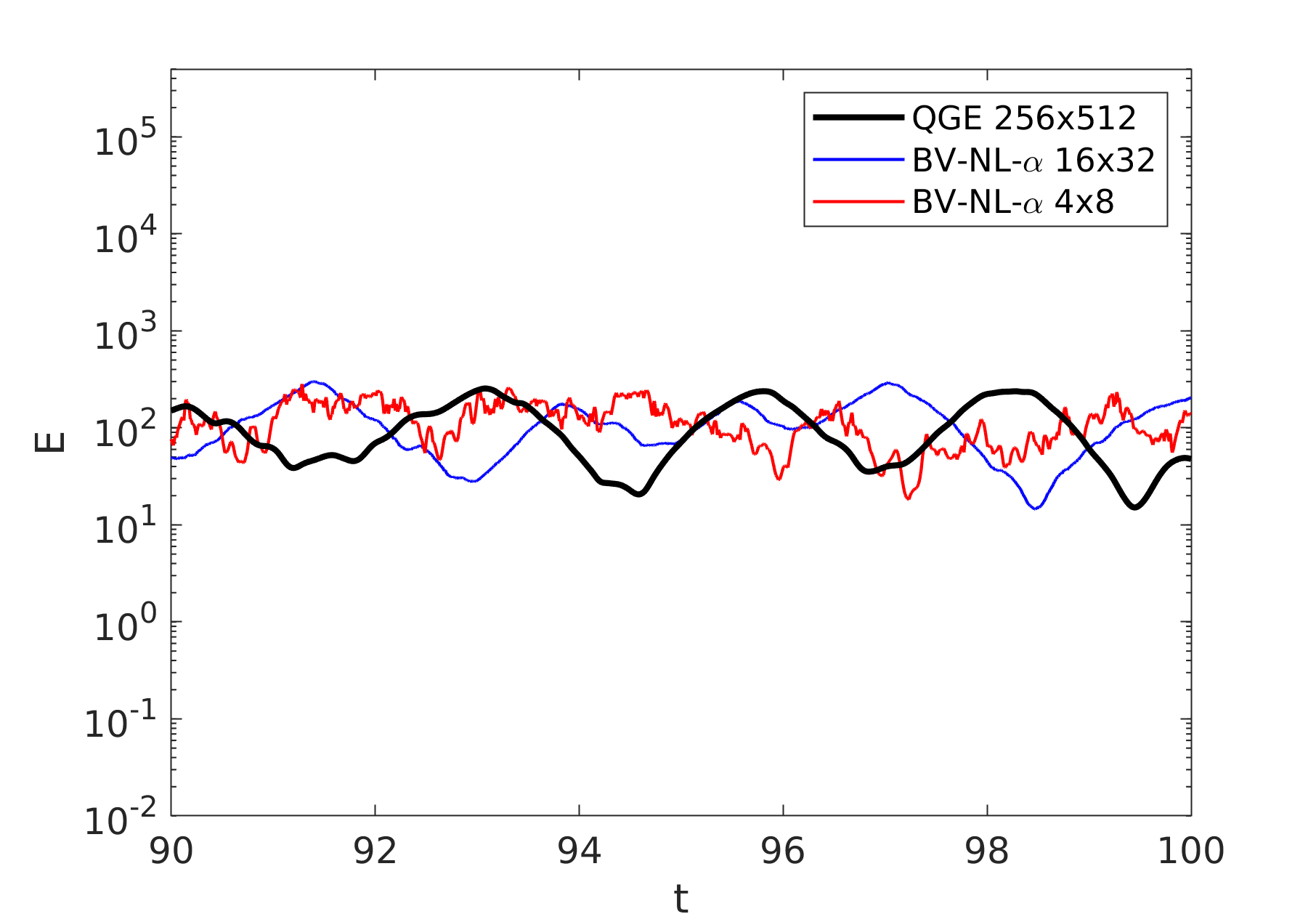}
      \end{overpic}
\caption{Case 2: time evolution of the kinetic energy \eqref{eq:kin_energy}
computed by the BV-NL-$\alpha$ model with meshes $16 \times 32$ and  $4 \times 8$ and compared to the
true kinetic energy (left) and corresponding zoomed-in view (right). 
} 
\label{fig:Ek_a_second_BVNLa}
\end{figure}

Fig.~\ref{fig:BV_second_filter} shows $\widetilde{\psi}$ and $\widetilde{q}$ computed with the BV-$\alpha$ model
and meshes $16 \times 32$ and  $4 \times 8$.
We observe no improvement over the QGE model with either mesh. 
This is confirmed by the computed kinetic energy reported in Fig.~\ref{fig:Ek_a_second_BVa}.
The solution computed with  the BV-NL-$\alpha$ model and mesh $16 \times 32$ is
shown in Fig.~\ref{fig:BV_NL_second_filter} (a) and (b). It 
does represent a significant improvement: 
the computed $\widetilde{\psi}$ and $\widetilde{q}$  are in very good agreement with the exact solution
both in terms of pattern and magnitude. 
The kinetic energy computed with  the BV-NL-$\alpha$ model and mesh $16 \times 32$ is also in good agreement
with the exact kinetic energy. In fact, Fig.~\ref{fig:Ek_a_second_BVNLa} shows that the amplitudes and frequencies
of oscillation are comparable, although the phases are off. The accurate solution
computed by the BV-NL-$\alpha$ model comes at a fraction of the computational time required
by the exact solution: the total CPU time required by the QGE simulation with mesh $256 \times 512$ is about 95 hours, 
while the BV-NL-$\alpha$ model with mesh $16 \times 32$ takes about 1 hour. 
Even with the coarsest mesh under consideration (mesh $4 \times 8$) the BV-NL-$\alpha$ model captures
the solution pattern better than the QGE and BV-$\alpha$ models, although the magnitude of $\widetilde{\psi}$ is reduced
by the nonlinear filter. The average kinetic energy computed with the BV-NL-$\alpha$ model and mesh
$4 \times 8$ is also in good agreement with the true average kinetic energy. See Fig.~\ref{fig:Ek_a_second_BVNLa}.

%

\begin{figure}[htb!]
\centering
\begin{subfigure}{0.193\textwidth}
         \centering
         \includegraphics[width=\textwidth]{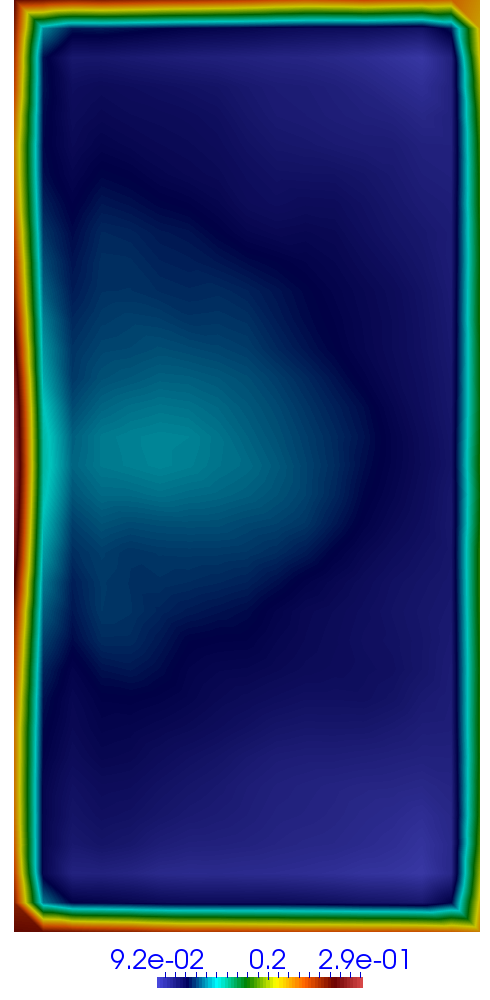}
         \caption{\scriptsize{$\widetilde{\psi}$, $16 \times 32$}}
     \end{subfigure}
     \begin{subfigure}{0.193\textwidth}
         \centering
         \includegraphics[width=\textwidth]{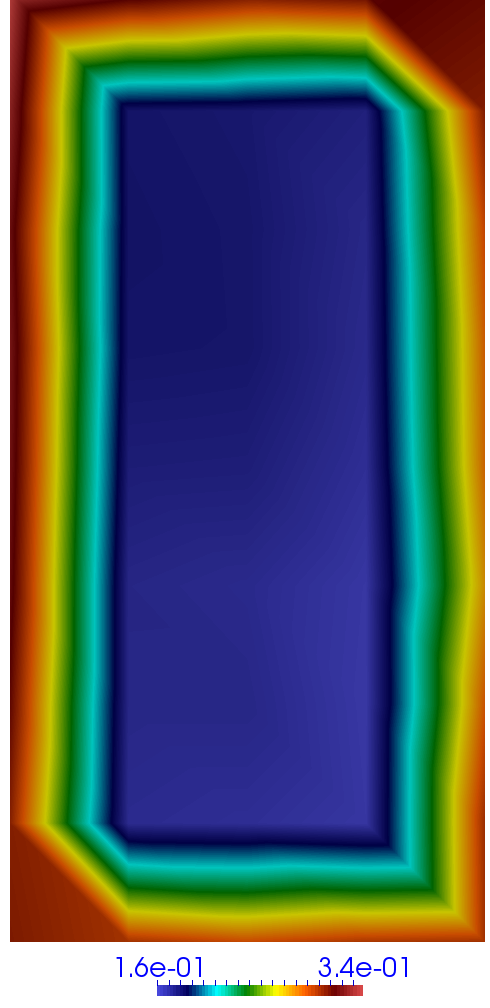}
         \caption{\scriptsize{$\widetilde{q}$, $4 \times 8$}}
     \end{subfigure}
\caption{Case 2: $\widetilde{a}$ computed with the BV-NL-$\alpha$ model and meshes $16 \times 32$ and  $4 \times 8$. }
\label{fig:indicator_func_second}
\end{figure}

Let us conclude with the visualization of the time-averaged indicator function $\widetilde{a}$ 
computed with the meshes $16 \times 32$ and $4 \times 8$ shown in Fig.~\ref{fig:indicator_func_second}.
We see that $\widetilde{a}$ takes its the largest values in one strip of cells close to the boundary, where
there are the largest gradients of the potential vorticity. Notice also that the maximum and minimum values 
of $\widetilde{a}$ become larger in absolute value when the mesh gets coarser, since more regularization is needed. 




\section{Conclusions and perspectives}\label{sec:conclusions}

We presented a nonlinear variant of the BV-$\alpha$ model, called BV-NL-$\alpha$, for the simulation of  
barotropic flows with under-refined meshes. 
To select the regions of the domain where filtering is needed, we employed 
a nonlinear differential low-pass filter. For the space discretization of the BV-NL-$\alpha$ model, 
we chose a Finite Volume method that has the advantage of enforcing conservation of quantities at the discrete level.

We showed the effectiveness of our approach through a computational study for the double-gyre wind forcing benchmark. 
We considered two different parameters setting: i) $Ro = 0.0036$ and $Re = 450$ and ii) $Ro = 0.008$ and $Re = 1000$. 
In both cases, when coarse meshes are considered the BV-NL-$\alpha$ model provides more accurate results than 
the QGE and BV-$\alpha$ models.

This work could be expanded in different directions. A sensitivity analysis for the filtering radius \cite{Bertagna2019}
would help us understand how to obtain the most accurate results when compared to a direct
numerical simulation. Moreover, it would be interesting to test the performance of a class of deconvolution-based 
indicator functions and to implement an efficient algorithm called Evolve-Filter-Relax, 
which proved to work well for the Leray-$\alpha$ model \cite{BQV, Girfoglio2019, Girfoglio2021a, Girfoglio2021b, Girfoglio2021c, Girfoglio_JCP, Strazzullo2021}. 
Thus, we believe they could be successful also for the BV-NL-$\alpha$ model. 



\section{Acknowledgements}\label{sec:acknowledgements}
We acknowledge the support provided by the European Research Council Executive Agency by the Consolidator Grant project AROMA-CFD ``Advanced Reduced Order Methods with Applications in Computational Fluid Dynamics" - GA 681447, H2020-ERC CoG 2015 AROMA-CFD, PI G. Rozza, and INdAM-GNCS 2019-2020 projects.
This work was also partially supported by US National Science Foundation through grant DMS-1953535. 
A.Quaini also acknowledges support from the Radcliffe Institute for Advanced Study at Harvard University where
she has been a 2021-2022 William and Flora Hewlett Foundation Fellow.

\bibliographystyle{plain}
\bibliography{QGE} 

\end{document}